\documentclass[12pt]{amsart}
\usepackage{graphicx} 
\usepackage{xcolor}
\usepackage{amssymb}
\usepackage{verbatim}
\usepackage[labelformat=simple]{subcaption}	
\usepackage[T1]{fontenc}
\usepackage[margin=1in]{geometry}
\usepackage{braket}
\usepackage{setspace}
\usepackage{amsmath,amsthm}
\usepackage{tikz}

\newcommand{\dee}{\ensuremath{\textrm{d}}}

\newcommand{\fdf}[1]{\ensuremath{ \frac{\dee}{\dee #1}}}
\newcommand{\fdfn}[2]{\ensuremath{ \frac{\dee^{#2}}{\dee #1^{#2}}}}

\newcommand{\fdn}[3]{\ensuremath{ \frac{\dee^{#3} #1}{ \dee #2^{#3}}}}
\newcommand{\inty}[4]{\ensuremath{ \int_{#1}^{#2} \! #3 \, \dee#4 }}

\newcommand{\ip}[2]{\ensuremath{ \left< \left. #1 \right| #2 \right> } }

\DeclareMathOperator{\DoS}{DoS}
\DeclareMathOperator{\Tr}{Tr}

\newtheorem{theorem}{Theorem}

\newtheorem{lemma}{Lemma}

\numberwithin{lemma}{section}
\numberwithin{example}{section}
\numberwithin{proposition}{section}
\numberwithin{equation}{section}
\numberwithin{theorem}{section}
\numberwithin{corollary}{section}
\numberwithin{remark}{section}
\numberwithin{definition}{section}
\numberwithin{assumption}{section}
\newtheorem*{definition*}{Definition}		
\newtheorem*{assumption*}{Assumption}
\newtheorem*{remark*}{Remark}
\newtheorem*{theorem*}{Theorem}
\newtheorem*{lemma*}{Lemma}
\newtheorem*{proposition*}{Proposition}
\newtheorem*{corollary*}{Corollary}
\newtheorem*{example*}{Example}
\newtheorem*{conjecture*}{Conjecture}

\title{Discrete-to-continuum convergence of the density of states for Mathieu's equation}

\author{Peter Hofhansel, Alexander B. Watson}
\date{\today}

\begin{document}

\begin{abstract}
The density of states of a self-adjoint operator generalizes the eigenvalue distribution of a Hermitian matrix. We prove convergence of the density of states for a tight-binding model with a slowly-varying periodic potential to the density of states of its continuum approximation, a Mathieu-type equation.
\end{abstract}

\maketitle

\tableofcontents

\section{Introduction} \label{sec:intro}

Discrete tight-binding models are used throughout condensed matter physics to model electrons in materials. When carefully parameterized by density functional theory (DFT) computations, tight-binding models are expected to retain excellent accuracy relative to more fundamental continuum Schr\"odinger models \cite{kaxiras_joannopoulos_2019,ashcroft_mermin}.

Continuum PDE approximations to discrete tight-binding models have attracted significant attention in recent years since they allow for further reduction in model complexity, an especially important consideration when parameterizing models of interacting electrons. Such PDE approximations play an essential role in predicting the electronic properties of moir\'e materials such as twisted bilayer graphene \cite{Bistritzer2011,Watson2022,Quinn2025,cances2023,cances2023semiclassical,Cances_secondorder}, which have attracted huge attention in recent years for their remarkable electronic properties \cite{Cao2018,Cao2018a,kim2021,ke2021}.

The goal of the present work is to provide a simple yet rigorous justification of the discrete-to-continuum approximation for the density of states: roughly speaking, the distribution of eigenvalues. For simplicity and clarity, here we consider simple model equations: the 1D discrete and continuum Mathieu equations, but we expect that the insights can be generalized to higher dimensions and to models of moir\'e materials. Convergence of the density of states to those of effective reduced models when the layers are incommensurate has been considered in \cite{Cances2017a,Massatt2017,MassattCarrLuskinOrtner2018,Massatt2020,Massatt2023,DAI2023112204,Wang2025}.

The structure of the work is as follows. We first review the limit for the discrete Laplacian, proving convergence of the density of states to that of the continuum Laplacian. We consider this case in Section \ref{sec:laplacian}. We then consider a discrete Mathieu equation: the discrete Laplacian with a periodic potential whose period is much larger than the lattice constant. Here we prove convergence of the density of states to that of the continuum Laplacian with a 1-periodic potential, i.e., a continuum Mathieu equation. We consider this case in Section \ref{sec:mathieu}. In each case, we introduce the model, provide a formal derivation of its continuum approximation, derive the density of states of the discrete model and its continuum approximation, and then state a theorem guaranteeing convergence of the density of states. Our theorems are supported by numerical computations comparing the discrete and continuum density of states distributions. All proofs are provided in the Appendices. 

\subsection*{Acknowledgements} This work was completed during PH's REU at the University of Minnesota during Summer 2025 supported by the Undergraduate Research Opportunies Program (UROP). ABW's research was supported in part by NSF grant DMS-2406981.

\section{From discrete Laplacian to continuum Laplacian} \label{sec:laplacian}

We start by studying the discrete-to-continuum limit for perhaps the simplest possible model: the 1D discrete Laplacian.

\subsection{Derivation of continuum limit via formal multiple-scales analysis} \label{sec:laplacian_formal}

We start by deriving the continuum limit by a systematic but nonrigorous calculation. The discrete Laplacian is the operator
\begin{equation} \label{eq:disc_lap_0}
    (H_d\psi)_n := -(\psi_{n+1} + \psi_{n-1} - 2\psi_{n})
\end{equation}
acting on $\psi \in \ell^2(\mathbb{Z})$. The continuum limit of this operator emerges when we consider the operator acting on functions which vary appreciably only over many grid points. We can express this precisely by assuming that
\begin{equation} \label{eq:ansatz_0}
    \psi_n = \left. \phi(x) \right|_{x = \epsilon n},
\end{equation}
where $\phi \in \mathcal{S}(\mathbb{R})$ (Schwartz class) is assumed to be a fixed smooth function, while $\epsilon = \frac{1}{N}$, where $N$ is the number of grid points over which the wave-function is assumed to vary appreciably. The continuum limit is then expected to emerge when $\epsilon \ll 1$, or equivalently, $N \rightarrow \infty$. Substituting \eqref{eq:ansatz_0} into \eqref{eq:disc_lap_0} we obtain
\begin{equation} 
    (H_d\psi)_n = -(\phi(\epsilon(n+1)) + \phi(\epsilon(n-1)) - 2\phi(\epsilon n)) = - \epsilon^2 \fdn{\phi}{x}{2}(\epsilon n) + O(\epsilon^4).
\end{equation}
We are thus motivated to introduce the continuum limit operator
\begin{equation} \label{eq:H_c definition}
    (H_c \phi)(x) := - \fdn{\phi}{x}{2}(x)
\end{equation}
acting on twice differentiable functions $\phi \in L^2(\mathbb{R})$, and to expect that eigenfunctions of $H_d$ with eigenvalue $E_d$ will converge to eigenfunctions of $H_c$ with eigenvalue $E_c := \frac{E_d}{\epsilon^2}$ in the limit $\epsilon \downarrow 0$.

\subsection{Density of states for discrete model} \label{sec:laplacian_discrete_DoS}

We now introduce the density of states (DoS) of the operator $H_d$ and derive a convenient formula for comparing this DoS with the DoS of the operator $H_c$, which we will introduce in Section \ref{sec:laplacian_continuum_DoS}. The DoS is a convenient proxy for the distribution of eigenvalues which is more convenient for numerical computation and analytical proofs, especially for infinite-dimensional operators whose spectrum may not be purely discrete. 

Recall that for Hermitian matrices $H$, we have the spectral theorem so that $H$ can be written 
\begin{equation*}
    H = U D U^{-1},
\end{equation*}
where $U$ is the unitary matrix whose columns are the normalized eigenvectors of $H$, and $D$ is the diagonal matrix whose entries are the real eigenvalues of $H$. We can then define functions of $H$ by
\begin{equation} \label{eq:f_H}
    f(H) := U f(D) U^{-1}
\end{equation}
whenever the right-hand side makes sense. The function
\begin{equation*}
    L : \mathbb{R} \rightarrow \mathbb{R}, \quad L(x) := (1 + x^2)^{-1} = \frac{-i}{2} \left( (x - i)^{-1} - (x + i)^{-1} \right)
\end{equation*}
is known as the Lorentzian. We can now define the DoS of an $N \times N$ Hermitian matrix $H$ to be the function 
\begin{equation} \label{eq:DOS}
    \DoS : \mathbb{R} \rightarrow \mathbb{R}, \quad \DoS(\mu) := \frac{1}{N} \Tr((1+(H - \mu )^2)^{-1},
\end{equation}
where $\Tr$ denotes the matrix trace. Note that the right-hand side makes sense since $H$ is Hermitian and hence has purely real eigenvalues.

We now consider the discrete Laplacian operator $H_d$ \eqref{eq:disc_lap_0} applied to complex vectors $\psi = (\psi_n)_{-\frac{N-1}{2} \leq n \leq \frac{N-1}{2}}$ with dimension $N$ where $N$ is odd, subject to periodic boundary conditions
$\psi_{\frac{N-1}{2}+1}=\psi_{-\frac{N-1}{2}}$ and $\psi_{-\frac{N-1}{2}-1}=\psi_{\frac{N-1}{2}}$. 
We note that the operator $H_d$ has matrix representation 
\begin{equation}
    H_d = \begin{bmatrix}
        2 & -1 & 0 & ... & ... & -1\\ -1 & 2 & -1 & 0 & ... & ...\\ 0 & -1 & 2 & -1 & 0 & ... \\ ... & ... & ... & ... & ... & ... \\ -1 & 0 & ... & ...& -1 & 2
    \end{bmatrix}.
\end{equation}
This is clearly a real symmetric matrix and is therefore Hermitian. 

The easiest basis to compute the trace of the Lorentzian and therefore the DoS of $H_d$ is the Fourier basis.
So we compute the Fourier Transform of $H_d$. For odd $N$ the Fourier modes have frequencies $\xi \in \set{-\frac{(N-1)\pi}{N},-\frac{(N-3)\pi}{N},...,\frac{(N-1)\pi}{N}}$. The Fourier Transform is
\begin{equation}
    (F\psi)_\xi= \frac{1}{\sqrt{N}}\sum_{n=-\frac{N-1}{2}}^{\frac{N-1}{2}}{\psi_ne^{-i\xi n}}, \quad (F^{-1}\hat{\psi})_n= \frac{1}{\sqrt{N}}\sum_{\xi}^{}{\hat{\psi}_\xi e^{i\xi n}}.
\end{equation}
So $H_d$ in the Fourier basis is 
\begin{multline}
    (D_d \hat{\psi})(\xi) := (FH_dF^{-1}\hat{\psi})_\xi = \frac{1}{N} \sum_{n=-\frac{N-1}{2}}^{\frac{N-1}{2}}\, e^{-i\xi n} \sum_{\xi'}^{}\, \left(-e^{i\xi'} 
    - e^{-i \xi'} + 2 \right)e^{i\xi' n}\hat{\psi}_{\xi'} \\ = \sum_{\xi'}^{}\, \hat{\psi}_{\xi'} \left(2- 2\cos \left( \xi' \right)\right) 
    \frac{1}{N} \sum_{n=-\frac{N-1}{2}}^{\frac{N-1}{2}} e^{i(\xi'-\xi)n}  = \hat{\psi}_{\xi}\left(2 - 2\cos\left(\xi\right)\right).
\end{multline}
In matrix form this is 
\begin{equation}
    FH_dF^{-1} = \begin{bmatrix}
... & ... & ... & ... & ... \\ ... & 2 \left( 1 - \cos\left(- \frac{2 \pi}{N} \right)\right) & 0  & ... & ... \\ ... & 0 & 0 & 0  & ... \\ ... & ... & 0  & 2 \left( 1 - \cos\left(\frac{2 \pi}{N} \right)\right) & ... \\ ... & ... & ... & ... & ...
\end{bmatrix}.
\end{equation}

Therefore applying the definition of the density of states \eqref{eq:DOS} and of the functional calculus \eqref{eq:f_H} we have
\begin{equation}
    \DoS_d(\mu) = \frac{1}{N} \Tr(1+(H_d - \mu )^2)^{-1} = \frac{1}{N} \Tr(F^{-1} \left(1+(F H_d F^{-1} - \mu )^2)^{-1} F\right), 
\end{equation}
then using $\Tr(AB) = \Tr(BA)$ we have
\begin{equation}
    = \frac{1}{N} \Tr(1+(F H_d F^{-1} - \mu )^2)^{-1} = \frac{1}{N} \sum_{n=-\frac{N-1}{2}}^{\frac{N-1}{2}}\frac{1}{1+(2 -2 \cos(\frac{2\pi n}{N}) - \mu )^2},
\end{equation}
which is nothing but a Riemann sum of a smooth function. Taking the limit $N \rightarrow \infty$ we obtain 
\begin{equation} \label{eq:disc_Lap_DoS}
    \DoS_d(\mu) = \frac{1}{2\pi} \int_{-\pi}^{\pi} \frac{1}{1+(2 -2 \cos(x) - \mu )^2} \text{d}x.
\end{equation}
We will show that this quantity, when appropriately scaled with $\epsilon$, converges to the analogous quantity for the continuum limit in Section \ref{sec:laplacian_convergence}. We derive the analogous quantity for the continuum limit in Section \ref{sec:laplacian_continuum_DoS}.

\subsection{Density of states for continuum model} \label{sec:laplacian_continuum_DoS}

We now define the DoS for operators without a finite grid and derive a convenient formula for the DoS of $H_c$ to compare to the DoS of $H_d$ previously calculated. We consider $H_c$ defined in (\ref{eq:H_c definition}) applied to functions on the interval $\left(-\frac{N}{2},\frac{N}{2}\right)$ with periodic boundary conditions. 
To do this as before we find the matrix form of $H_c$ in the Fourier basis.

We consider the one-dimensional continuum Laplacian operator \eqref{eq:H_c definition}
\begin{equation}
    ( H_c \phi )(x) := - \fdn{\phi}{x}{2}(x)
\end{equation}
with the domain 
\begin{equation}
    D(H_c) := \left\{ f \in L^2\left(I\right) : f', f'' \in L^2(I), f\left(-\frac{N}{2}\right) = f\left(\frac{N}{2}\right), f'\left(-\frac{N}{2}\right) = f'\left(\frac{N}{2}\right) \right\},
\end{equation}
where $I :=  \left(-\frac{N}{2},\frac{N}{2}\right)$. This operator is self-adjoint, with an orthonormal basis of $L^2(I)$ consisting of eigenfunctions of $H_c$ provided by the Fourier basis
\begin{equation} \label{eq:xi_Fourier}
    \left\{ \frac{1}{\sqrt{N}} e^{i \xi x} : \xi = \frac{2 \pi k}{N}, k \in \mathbb{Z} \right\}.
\end{equation}
Equivalently, $H_c$ is diagonalized by the Fourier transform
\begin{equation}
    (F\phi)(\xi) = \frac{1}{N} \inty{-\frac{N}{2}}{\frac{N}{2}}{e^{-i\xi x} \phi(x)}{x}, \quad (F^{-1}\hat{\phi})(x) = \sum_{\xi}^{}\,e^{i\xi x} \hat{\phi}(\xi),
\end{equation}
where $\xi = \frac{2 \pi k}{N}, k \in \mathbb{Z}$, since
\begin{equation}
    \begin{split}
        \left( D_c \hat{\phi} \right)(\xi) := (FH_cF^{-1}\hat{\phi})(\xi) =& \frac{1}{N} \inty{-\frac{N}{2}}{\frac{N}{2}}{e^{-i\xi x} \xi'^2 \sum_{\xi'}^{}\,e^{i\xi' x}\hat{\phi}(\xi') }{x} \\ =& \sum_{\xi'}^{}\,\hat{\phi}(\xi') \xi'^2 \frac{1}{N} \inty{-\frac{N}{2}}{\frac{N}{2}}{ e^{i(\xi' - \xi)x} }{x} = \hat{\phi}(\xi) \xi^2.
    \end{split}
\end{equation}

We can now define functions of $H_c$ as in the discrete case by
\begin{equation}
    f(H_c) := F^{-1} f\left(D_c\right) F,
\end{equation}
and their trace (whenever the sum converges) by
\begin{equation}
    \Tr( f(H_c) ) := \sum_{\xi} \ip{\frac{1}{\sqrt{N}} e^{i \xi x}}{f(H_c) \frac{1}{\sqrt{N}} e^{i \xi x}} = \sum_{\xi} f( D_c(\xi) ) = \sum_{\xi} f( \xi^2 ).
\end{equation}

We now define the DoS of a continuous operator to match as much as possible the definition for a discrete operator. We use the definition of the DoS for a discrete operator with the new definition of the trace and we replace the number of grid-points in the formula with the length of the domain of the functions in the domain of the operator, so that
\begin{multline}
    \DoS_c(\mu) := \frac{1}{N}\Tr((1+(H_c - \mu )^2)^{-1} = \frac{1}{N} \sum_{\xi} ((1+(\xi^2 - \mu )^2)^{-1} = \frac{1}{N} \sum_{k \in \mathbb{Z}} \frac{1}{1+\left(\left( \frac{2\pi k}{N}\right)^2 - \mu \right)^2},
\end{multline}
which clearly converges absolutely.

Since this is a Riemann sum of a smooth function, we can take the limit $N \rightarrow \infty$ to obtain 
\begin{equation} \label{eq:continuum diagonal Dos}
    \DoS_c(\mu) = \frac{1}{2\pi} \int_{-\infty}^{\infty}\,\frac{1}{1+(k^2 - \mu )^2} \text{d}k.
\end{equation}

\subsection{Convergence theorem} \label{sec:laplacian_convergence}

\begin{figure}[!htb]
\centering
\includegraphics[width=\columnwidth]{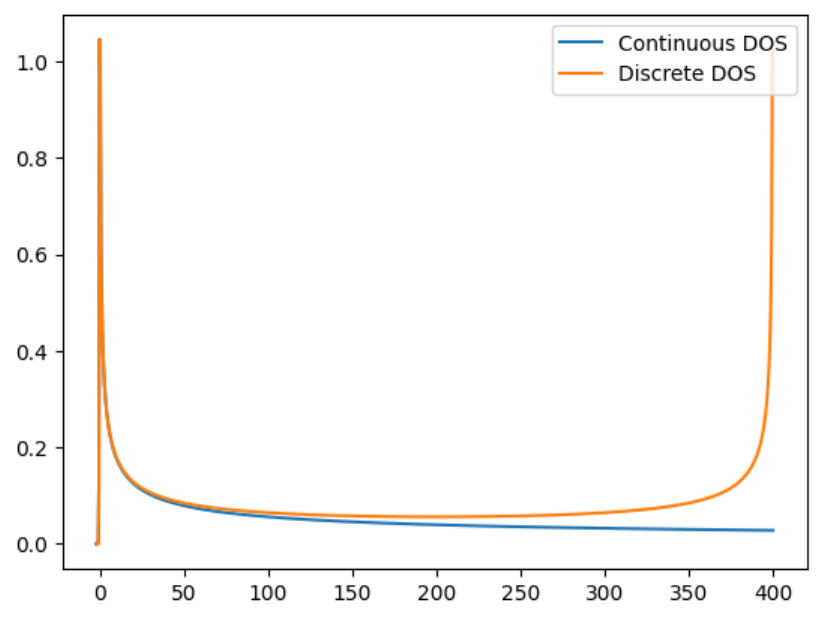}
    \caption{\label{fig:DiagonalCase} Numerically computed re-scaled density of states \eqref{eq:Lap_disc_DoS} for the discrete Laplacian \eqref{eq:disc_lap_0}, plotted as a function of $\mu$, compared with the numerically computed density of states \eqref{eq:continuum diagonal Dos} of the continuum Laplacian \eqref{eq:H_c definition}, showing the expected agreement near the bottom of the spectrum.}
\end{figure}

Near the bottom of the spectrum, the discrete DoS \eqref{eq:disc_Lap_DoS} can be accurately approximated by the continuum DoS \eqref{eq:continuum diagonal Dos}. To make this precise, we consider the $\epsilon$-scaled discrete DoS
\begin{equation} \label{eq:Lap_disc_DoS}
    \text{DoS}_{d}^{\epsilon}(\mu) := \frac{1}{\epsilon} \lim_{N \rightarrow \infty} \frac{1}{N} \Tr \left( 1 + \left(\frac{H_d}{\epsilon^2} - \mu \right)^2 \right)^{-1} = \frac{1}{2\pi\epsilon} \int_{-\pi}^{\pi}\,\frac{1}{1+\left(\frac{2 - 2 \cos(k)}{\epsilon^2} - \mu \right)^2} \text{d}k.
\end{equation}
Our goal now is to prove the following theorem:
\begin{theorem} \label{thm:Laplace_theorem}
For each fixed $\mu \in \mathbb{R}$, there exist positive constants $C > 0$ and $\epsilon_0 > 0$ such that, for all $0 < \epsilon < \epsilon_0$, we have
\begin{equation}
    \left| \text{\emph{DoS}}_d^\epsilon(\mu) - \text{\emph{DoS}}_c(\mu) \right| \leq C \epsilon^{\frac{1}{11}}.
\end{equation}
In particular,
\begin{equation}
    \lim_{\epsilon \rightarrow 0} \text{\emph{DoS}}_d^\epsilon(\mu) = \text{\emph{DoS}}_c(\mu).
\end{equation}
\end{theorem}
To see that this is plausible, note that if we change variables we have 
\begin{equation}
    \text{DoS}_{d}^{\epsilon}(\mu) = \frac{1}{2\pi} \int_{-\frac{\pi}{\epsilon}}^{\frac{\pi}{\epsilon}}\,\frac{1}{1+\left(\frac{2 - 2 \cos(\epsilon k)}{\epsilon^2} - \mu \right)^2} \text{d}k,
\end{equation}
where the limits converge to $\pm \infty$, and Taylor-expansion shows that $\frac{2 - 2 \cos(\epsilon k)}{\epsilon^2} \approx k^2 + O(\epsilon^2 k^4)$. However, this calculation does not constitute a full proof. We provide a detailed proof in Appendix \ref{sec:laplacian_proofs}. Numerical results for the density of states calculated using both formulae are shown in Figure \ref{fig:DiagonalCase}.

\section{From discrete to continuum Mathieu equation} \label{sec:mathieu}

We now consider the discrete-to-continuum limit for the discrete Mathieu equation.

\subsection{Derivation of continuum limit via formal multiple-scales analysis} \label{sec:mathieu_forma}

We again start by deriving the continuum limit by a systematic but non-rigorous calculation. The discrete Mathieu operator is 
\begin{equation}\label{eq:disc_matheiu}
    (H_d\psi)_n = -(\psi_{n+1} + \psi_{n-1} - 2\psi_{n}) + 2\lambda\epsilon^2\cos(2\pi\epsilon n)\psi_n
\end{equation}
acting on $\psi \in \ell^2(\mathbb{Z})$, where $\lambda$ is a real constant. This operator is periodic with period $\frac{1}{\epsilon}$ when $\frac{1}{\epsilon}$ is an integer and the continuum limit of this operator emerges when there are many grid points in a single period of the operator and the operator acts on functions which only vary appreciably on the same scale as the period of $H_d$.  We can express this precisely by assuming that
\begin{equation} \label{eq:ansatz}
    \psi_n = \left. \phi(x) \right|_{x = \epsilon n},
\end{equation}
where $\phi \in \mathcal{S}(\mathbb{R})$ is assumed to be a fixed smooth function. 
The continuum limit is then expected to emerge when $\epsilon \ll 1$. 
Substituting \eqref{eq:ansatz} into \eqref{eq:disc_matheiu} we obtain
\begin{multline} 
    (H_d\psi)_n = -(\phi(\epsilon(n+1)) + \phi(\epsilon(n-1)) - 2\phi(\epsilon n)) + 2\lambda\epsilon^2\cos(2\pi\epsilon n)\phi(\epsilon n) \\ = - \epsilon^2 \fdn{\phi}{x}{2}(\epsilon n) + 2\lambda\epsilon^2\cos(2\pi\epsilon n)\phi(\epsilon n) + O(\epsilon^4).
\end{multline}
We are thus motivated to introduce the continuum limit operator
\begin{equation} \label{eq:H_c_Mathieu}
    (H_c \phi)(x) := - \fdn{\phi}{x}{2}(x) + 2\lambda \cos(2\pi x)\phi(x)
\end{equation}
acting on $\phi \in L^2(\mathbb{R})$, and expect that eigenfunctions of $H_d$ with eigenvalue $E_d$ will converge to eigenfunctions of $H_c$ with eigenvalue $E_c := \frac{E_d}{\epsilon^2}$ in the limit $\epsilon \downarrow 0$.

\subsection{Density of states for discrete model} \label{sec:mathieu_discrete_DoS}

In this section, we derive a convenient formula for the DoS of $H_d$ to compare with the DoS of the operator $H_c$, which we will find in Section \ref{sec:mathieu_continuum_DoS}. 
We use the same definition of the density of states as when studying the Laplacian 
\begin{equation}
    \DoS(\mu) = \frac{1}{m} \Tr((1+(H - \mu )^2)^{-1}
\end{equation}
for an $m \times m$ Hermitian $H$.
In this case as we are dealing with a periodic operator we study the operator over $M$ periods and therefore apply it to vectors of size $MN$ where $N = \frac{1}{\epsilon}$ is the period of $H_d$.
The boundary conditions are $\psi_{MN}=\psi_{0}$ and $\psi_{MN-1}=\psi_{-1}$. 

We now consider the translation operator $T$ defined by
\begin{equation}
    (T\psi)_n = \psi_{n+N}, \quad (T^{-1}\psi)_n = \psi_{n-N}.
\end{equation}
The matrix representation of $T$ is
\begin{equation}
    T_{mn} = \begin{cases}
1, & \text{if } m = n+N\\
0, & \text{otherwise}
\end{cases}, \quad T^{-1}_{mn} = \begin{cases}
1, & \text{if } m = n-N\\
0, & \text{otherwise}
\end{cases}.
\end{equation}
From this we can see that $T^{-1}$ is the conjugate transpose of $T$. Therefore $T$ is unitary. It is easy to confirm that $H_d$ and $T$ commute
\begin{equation}
    (H_dT\psi)_n = -(\psi_{n+1+N}-2\psi_{n+N} + \psi_{n-1+N}) + 2\lambda\epsilon^2\cos(2\pi\epsilon n)\psi_{n+N}= (TH_d\psi)_n.
\end{equation}

The operator $H_d$ only differs from how it was defined in the previous section on the main diagonal and so is also Hermitian.
Since $T$ and $H_d$ are both normal and they commute, there exists a basis of simultaneous eigenvectors of $T$ and $H_d$. Therefore when seeking eigenvectors of $H_d$, we can restrict attention WLOG to eigenvectors of $T$. Since $T$ is unitary, its eigenvalues must have magnitude 1. We can therefore assume WLOG that eigenvectors of $T$ satisfy
\begin{equation}
    \psi_{n+N} = e^{ik}\psi_n
\end{equation}
for some $- \pi \leq k < \pi$. Since $(H_d\psi)_{MN}=(H_d\psi)_{0}$ considering only odd $M$
\begin{equation} \label{eq:allowed_k}
    \psi_{MN} = e^{iMk}\psi_{MN} \implies e^{iMk} =1 \implies k = \frac{2j\pi}{M},
\end{equation}
where $j$ is an integer satisfying
\begin{equation}
    \frac{-(M-1)}{2}\leq j \leq \frac{M-1}{2}.
\end{equation}
We can now write all eigenvectors of $T$ as 
\begin{equation}
    \psi_n = e^{\frac{ikn}{N}} p_n,
\end{equation}
where $p$ is periodic with period $N$.

We now introduce the Bloch transform operator, which converts the operator $H_d$ acting on $MN$-vectors into a direct sum of operators $H_d(k)$ acting on periodic $N$-vectors. We define the Bloch transform at wave-number $k$ by 
\begin{equation}\label{eq:discrete_bloch_def}
    \hat{\psi}_{k,n} = (S \psi)_{k,n} := \frac{1}{\sqrt{M}}\sum_{m=0}^{M-1}\,e^{- ik(m+\frac{n}{N})} \psi_{n+Nm},
\end{equation}
where $k$ ranges over the allowed values \eqref{eq:allowed_k}.
To see that $S$ maps to $N$-periodic vectors, note that for any $k$ we have
\begin{multline}
    (S \psi)_{k,n+N} = \sum_{m=1}^{M}\,e^{ik \left(m+1+\frac{n}{N} \right)} \psi_{n+N(m+1)} = \sum_{m=2}^{M}\,e^{ik\left(m+\frac{n}{N}\right)} \psi_{n+Nm} + e^{ik(M+1+\frac{n}{N})}\psi_{n+(M+1)N} = \\\sum_{m=2}^{M}\,e^{ik\left(m+\frac{n}{N}\right)} \psi_{n+Nm} + e^{ik\left(1+\frac{n}{N}\right)}\psi_{n+N} = \sum_{m=1}^{M}\,e^{ik\left(m+\frac{n}{N}\right)} \psi_{n+Nm} = (S\psi)_{k,n}.
\end{multline}
We therefore restrict $n$ to the set $\{0,...,N-1\}$ for Bloch-transformed vectors. Since there are $M$ allowed $k$ values, we get that $S : \mathbb{C}^{MN} \rightarrow \mathbb{C}^M \otimes \mathbb{C}^N$. The inverse map is defined by
\begin{equation}
    \psi_n = (S^{-1} \hat{\psi})_{n} := \frac{1}{\sqrt{M}}\sum_{k}\,e^{ik \frac{n}{N}} \hat{\psi}_{k,n},
\end{equation}
where $k$ is summed over the allowed values \eqref{eq:allowed_k}. The Bloch transform is unitary, and conjugates $H_d$ to a block-diagonal matrix with respect to $k$. To see this, first note that
\begin{equation} \label{eq:commute_T_0}
    H_d T = T H_d \iff T^{-1} H_d T = H_d \implies \left( H_d \right)_{n + N m, n' + N m} = \left( H_d \right)_{n, n'}
\end{equation}
for any integer $m$. We now have 
\begin{equation}
    ( S H_d S^{-1} \hat{\psi} )_{k,n} = \frac{1}{M} \sum_m e^{- i k \left(m + \frac{n}{N}\right)} \sum_{n'} H_{n + N m, n'} \sum_{k'} e^{i k' \frac{n'}{N}} \hat{\psi}_{k',n'}.
\end{equation}
Re-ordering the sums, we have
\begin{equation}
    ( S H_d S^{-1} \hat{\psi} )_{k,n} = \sum_{k'} \sum_{n'} \left( \hat{H}_d \right)_{k,n,k',n'} \hat{\psi}_{k',n'},
\end{equation}
where
\begin{equation}
    \left( \hat{H}_d \right)_{k,n,k',n'} := \frac{1}{M} \sum_m e^{- i k \left(m + \frac{n}{N}\right)} \left( H_d \right)_{n + N m, n'} e^{i k' \frac{n'}{N}}
\end{equation}
are the matrix elements of $H_d$ with respect to the Bloch-transformed basis. Using periodicity of $\hat{\psi}_{k',n'}$ in $n'$, we can re-write the right-hand side as
\begin{equation}
    ( S H_d S^{-1} \hat{\psi} )_{k,n} = \sum_{k'} \sum_{n'} \frac{1}{M} \sum_m e^{- i k \left(m + \frac{n}{N}\right)} \left( H_d \right)_{n + N m, n' + N m} e^{i k' \left(m + \frac{n'}{N}\right)} \hat{\psi}_{k',n'}.
\end{equation}
Using \eqref{eq:commute_T_0} we can then simplify as 
\begin{equation}
    = \sum_{n'} \sum_{k'} \frac{1}{M} \sum_m e^{- i k \left(m + \frac{n}{N}\right)} \left( H_d \right)_{n, n'} e^{i k' \left(m + \frac{n'}{N}\right)} \hat{\psi}_{k',n'},
\end{equation}
so that $m$ now only appears in the exponentials. Performing the sum over $m$ we obtain 
\begin{equation}
    = \sum_{k'} \sum_{n'} \delta_{k,k'} e^{- i k \frac{n}{N}} \left( H_d \right)_{n, n'} e^{i k' \frac{n'}{N}} \hat{\psi}_{k',n'} = \sum_{n'} e^{- i k \frac{n - n'}{N}} \left( H_d \right)_{n,n'} \hat{\psi}_{k,n'}.
\end{equation}
The fact that the right-hand side does not involve a sum over $k$ shows that with respect to the Bloch-transformed basis $H_d$ becomes block-diagonal with respect to $k$. The $N \times N$ Hamiltonians within each block are known as the Bloch Hamiltonians, and explicitly they act on periodic $N$-vectors $p$ by
\begin{equation}\label{eq:matheiu H_d(k)}
    \left( H_d(k) p \right)_n = -\left(e^{i\frac{k}{N}} p_{n+1} + e^{-i\frac{k}{N}}p_{n-1} - 2p_{n}\right)+2\lambda\epsilon^2\cos(2\pi\epsilon n)p_n.
\end{equation}

Since $H_d$ is block diagonal in this basis its trace is the sum of the traces of its blocks. Additionally the elementary matrix operations in the Lorentzian can be applied to each block individually and the result will still be block diagonal. So we now only need to compute the Lorentzian of each individual block. To do this we take the Fourier transform of the blocks.

We now take the Fourier transform of $H_d(k)$. For odd $N$ the Fourier modes have frequencies $\xi \in \set{-\frac{(N-1)\pi}{N},-\frac{(N-3)\pi}{N},...,\frac{(N-1)\pi}{N}}$. The Fourier transform is
\begin{equation}
    (F\psi)_\xi= \frac{1}{\sqrt{N}}\sum_{n=1}^{N}{\psi_ne^{-i\xi n}}, \quad (F^{-1}\hat{\psi})_n= \frac{1}{\sqrt{N}}\sum_{\xi}^{}{\hat{\psi}_\xi e^{i\xi n}}.
\end{equation}
So $H_d(k)$ in the Fourier basis is 
\begin{multline}
    (FH_d(k)F^{-1}p)_\xi = \frac{1}{N} \sum_{n=1}^{N}\, e^{-i\xi n} \sum_{\xi'}^{}\, \left(-e^{i\left(\frac{k}{N} + \xi'\right)} - e^{i\left(\frac{-k}{N} - \xi'\right)} + 2 + 2\lambda\epsilon^2\cos(2\pi\epsilon n)\right)e^{i\xi' n}\hat{p}_{\xi'} \\ = \sum_{\xi}^{}\, p_{\xi'} \left(2- 2\cos\left(\xi' + \frac{k}{N}\right)\right) \frac{1}{N} \sum_{n=1}^{N} e^{i(\xi'-\xi)n} + p_{\xi'} \lambda \epsilon^2 \frac{1}{N} \sum_{n=1}^{N} \left(e^{i(\frac{2\pi}{N} + \xi'-\xi)n} + e^{i(\frac{-2\pi}{N} + \xi'-\xi)n}\right) \\ = p_{\xi}\left(2 - 2\cos\left(\xi' + \frac{k}{N}\right)\right) + \lambda \epsilon^2 \left(p_{\xi + \frac{2\pi}{N}} + p_{\xi - \frac{2\pi}{N}}\right) =: \left( D_d(k) p \right)_\xi.
\end{multline}
In matrix form this is 
\begin{equation}\label{eq:matrix FH_dF_-1}
\left( D_d(k) p \right)_\xi = \begin{bmatrix}
... & ... & ... & ... & ... \\ ... & 2 \left( 1 - \cos\left( \frac{k}{N} - \frac{2 \pi}{N} \right)\right) & \epsilon^2 \lambda  & ... & ... \\ ... & \epsilon^2 \lambda & 2 \left( 1 - \cos\left( \frac{k}{N} \right)\right) & \epsilon^2 \lambda  & ... \\ ... & ... & \epsilon^2 \lambda  & 2 \left( 1 - \cos\left( \frac{k}{N} + \frac{2 \pi}{N} \right)\right) & ... \\ ... & ... & ... & ... & ...
\end{bmatrix}.
\end{equation}

Since $H_d$ is block diagonal we can write the trace of the Lorentzian of $H_d$ as the sum of the traces of the Lorentzian applied to $D_d(k)$ for each $k$ so that
\begin{equation}
    \DoS_d(\mu) := \frac{1}{NM} \Tr\left( 1 + (H_d - \mu)^2 \right) = \frac{1}{NM} \sum_{k}^{}\,\Tr(1+(D_d(k) - \mu )^2)^{-1}.
\end{equation}
We are studying the system in the limit as $N$ and $M$ approach $\infty$. Since $k$ has $M$ values evenly distributed between $-\pi$ and $\pi$, in the limit as $M$ approaches $\infty$,
\begin{equation}\label{eq:disc_Mathieu_DoS}
    \DoS_d(\mu) = \frac{1}{2\pi N} \int_{-\pi}^{\pi}\,\Tr(1+(D_d(k) - \mu )^2)^{-1} \text{d}k.
\end{equation}
We will show that this quantity, when appropriately scaled with $\epsilon$, converges to the analogous quantity for the continuum limit in Section \ref{sec:mathieu_convergence}.

\subsection{Density of states for continuum model} \label{sec:mathieu_continuum_DoS}

In this section we derive a convenient formula for the DoS of $H_c$ \eqref{eq:H_c_Mathieu} to compare to the DoS of $H_d$ previously calculated. In the continuum case we consider $H_c$ applied to functions on the interval $\left(-\frac{M}{2},\frac{M}{2}\right)$ with periodic boundary conditions. This interval corresponds to the discrete case with $NM$ grid points according to the relation between $\phi$ and $\psi$. 

More precisely, we consider the one-dimensional continuum Mathieu operator \eqref{eq:H_c_Mathieu}
\begin{equation}
    ( H_c \phi )(x) := - \fdn{\phi}{x}{2}(x) + 2 \lambda \cos(2 \pi x) \phi(x).
\end{equation}
Since the potential is bounded, this operator is self-adjoint with the domain 
\begin{equation}
    D(H_c) := \left\{ f \in L^2\left(I\right) : f', f'' \in L^2(I), f\left(-\frac{M}{2}\right) = f\left(\frac{M}{2}\right), f'\left(-\frac{M}{2}\right) = f'\left(\frac{M}{2}\right) \right\},
\end{equation}
where $I :=  \left(-\frac{M}{2},\frac{M}{2}\right)$. Unlike the continuum Laplacian, we do not have access to an explicit basis of eigenfunctions of $H_c$. However, we will see by a change of basis it is block-diagonal. We will then obtain a formula for the density of states as the sum of the traces of the Lorentzian applied to each block.  

Note first that $H_c$ commutes with the unitary operator $T$
\begin{equation}
    (T\phi)(x) = \phi(x+1),
\end{equation}
so it is natural to try to simultaneously diagonalize $H_c$ and $T$. Since $T$ is unitary, its eigenvalues must have magnitude 1. We can therefore assume WLOG that eigenfunctions of $T$ satisfy
\begin{equation}
    \phi(x+1) = e^{ik}\phi(x).
\end{equation}
Since $(H_c\phi)(\frac{M}{2})=(H_c\phi)(-\frac{M}{2})$ considering only odd $M$
\begin{equation}\label{eq:continuous allowed k}
    \phi\left(\frac{M}{2}\right)= e^{iMk}\phi\left(\frac{-M}{2}\right) \implies e^{iMk} =1 \implies k = \frac{2j\pi}{M} 
\end{equation}
\begin{equation}
    \frac{-(M-1)}{2}\leq j \leq \frac{M-1}{2}.
\end{equation}
We can write all eigenfunctions of $T$ as 
\begin{equation}
    \phi(x) = e^{ikx} p(x)
\end{equation}
for $p$ with period 1.

We now introduce the Bloch transform operator, which converts the operator $H_c$ acting on $M$-periodic functions into a direct sum of operators $H_c(k)$ acting on 1-periodic functions. We define the Bloch transform at wave-number $k$ by 
\begin{equation}\label{eq:continuous_bloch_def}
    (S\phi)_k(x) = \frac{1}{\sqrt{M}}\sum_{m=-\frac{M-1}{2}}^{\frac{M-1}{2}}\,e^{- ik(m+x)} \phi(x+m),
\end{equation}
where $k$ ranges over the allowed values \eqref{eq:continuous allowed k}.
To see that for periodic $\psi$ with period $M$ $S\psi$ has period 1 note that
\begin{multline}
    (S\phi)_k(x+1) = \frac{1}{\sqrt{M}}\sum_{m=-\frac{M-1}{2}}^{\frac{M-1}{2}}\,e^{- ik(m+x+1)} \phi(x+1+m) = \frac{1}{\sqrt{M}}\sum_{m=-\frac{M-3}{2}}^{\frac{M+1}{2}}\,e^{- ik(m+x)} \phi(x+m) = \\ \frac{1}{\sqrt{M}}\sum_{m=-\frac{M-1}{2}}^{\frac{M-3}{2}}\,e^{- ik(m+x)} \phi(x+m) = (S\phi)_k(x).
\end{multline}
Therefore $S\phi$ is a set of $M$ functions of period 1. The original function $\phi$ can be recovered from the set $S\phi = \hat{\phi}$ with $S^{-1}$ which is

\begin{equation}
    (S^{-1}\hat{\phi})(x) = \frac{1}{\sqrt{M}} \sum_{k}^{} e^{ikx} \hat{\phi}_k(x).
\end{equation}
In this basis $H_c$ is block-diagonal with the form 
\begin{equation}
    \begin{split}
        (SH_cS^{-1}\hat{\phi})_k(x) &= \frac{1}{M} \sum_{m=-\frac{M-1}{2}}^{\frac{M-1}{2}}\,e^{- ik(m+x)} \left.\left(-\fdfn{y}{2} + 2\lambda\cos(2\pi y)\right) \left(\sum_{k'}^{} e^{ik'y} \hat{\phi}_{k'}(y)\right) \right|_{y = x + m} \\ &= \frac{1}{M} \sum_{k'}^{}\,e^{ix(k'-k)} \sum_{m=-\frac{M-1}{2}}^{\frac{M-1}{2}}\, e^{i(k'-k)m} \left(-\left(\fdf{x}+ik\right)^2 + 2\lambda\cos(2\pi (x+m)) \right)\hat{\phi}_{k'}(x+m) \\ &= \sum_{k'}^{}\,e^{ix(k'-k)} \left(-\left(\fdf{x}+ik\right)^2 + 2\lambda\cos(2\pi x) \right) \hat{\phi}_{k'}(x)\frac{1}{M} \sum_{m=-\frac{M-1}{2}}^{\frac{M-1}{2}}\, e^{i(k'-k)m} \\ &= \sum_{k'}^{}\,e^{ix(k'-k)} \left(-\left(\fdf{x}+ik\right)^2 + 2\lambda\cos(2\pi x) \right) \hat{\phi}_{k'}(x)\delta_{k,k'} \\ &= \left(-\left(\fdf{x}+ik\right)^2 + 2\lambda\cos(2\pi x) \right) \hat{\phi}_{k'}(x).
    \end{split}
\end{equation}
The fact that the right-hand side does not involve a sum over $k$ shows that with respect to the Bloch-transformed basis $H_c$ becomes block-diagonal with respect to $k$ and the block of the operator corresponding to $k$ is 
\begin{equation}\label{eq:mathieu H_c(k)}
    H_c(k) = -\left(\fdf{x}+ik\right)^2 + 2\lambda\cos(2\pi x).
\end{equation}

Since $H_c$ is block diagonal in this basis we can compute functions of $H_c$ block-wise and compute the trace by summing the traces of each block. We now compute the matrix form of each bloch in the Fourier basis. 
To do this we take the Fourier transform of the blocks using the Fourier transformation definition
\begin{equation}
    (F\phi)(\xi) = \int_{0}^{1}\,e^{-i\xi x} \phi(x) \text{d}x, \quad (F^{-1}\hat{\phi})(x) = \sum_{\xi}^{}\,e^{i\xi x} \hat{\phi}(\xi),
\end{equation}
where the $\xi$s are as in \eqref{eq:xi_Fourier}. The result is
\begin{multline}
    (FH_c(k)F^{-1}p)(\xi) = \int_{0}^{1}\,e^{-i\xi x}(\xi'^2 + k^2 + 2k\xi' + 2\lambda\cos(2\pi x)) \sum_{\xi'}^{}\,e^{i\xi' x}\hat{p}(\xi') \text{d}x \\ =  \sum_{\xi'}^{}\,\hat{p}(\xi')(k+\xi')^2 \int_{0}^{1}\, e^{i(\xi' - \xi)x} + \hat{p}(\xi')\lambda \int_{0}^{1}\, e^{i(2\pi + \xi' - \xi)x} + e^{i(-2\pi + \xi' - \xi)x}\text{d}x \\ = \hat{p}(\xi)(k+\xi)^2 + \lambda \hat{p}(\xi + 2\pi) + \lambda \hat{p}(\xi - 2\pi) =: (D_c(k) p)(\xi),
\end{multline}
which in matrix form is 
\begin{equation}\label{eq:matrix FH_cF_-1}
    D_c(k) = \begin{bmatrix}
... & ... & ... & ... & ... \\ ... & (k - 2 \pi)^2 & \lambda  & ... & ... \\ ... & \lambda & k^2 & \lambda  & ... \\ ... & ... & \lambda  & (k + 2 \pi)^2 & ... \\ ... & ... & ... & ... & ...
\end{bmatrix}.
\end{equation}

Applying this to the definition of the DoS we obtain 
\begin{equation}
    \DoS_c(\mu) := \frac{1}{M} \Tr((1+(H_c - \mu )^2)^{-1} = \frac{1}{M} \sum_{k}^{}\,\Tr(1+(D_c(k) - \mu )^2)^{-1},
\end{equation}
where the $k$ values are as in \eqref{eq:continuous allowed k}. Since $k$ has values evenly distributed between $-\pi$ and $\pi$ and the trace is clearly smooth as a function of $k$, in the limit as $M$ approaches $\infty$ we have
\begin{equation}\label{eq:cont_Mathieu_DoS}
    \DoS_c(\mu) = \frac{1}{2\pi} \int_{-\pi}^{\pi}\,\Tr(1+(D_c(k) - \mu )^2)^{-1} \text{d}k.
\end{equation}

\subsection{Convergence theorem} \label{sec:mathieu_convergence}

Near the bottom of the spectrum, we expect that the discrete DoS \eqref{eq:disc_Mathieu_DoS} can be accurately approximated by the continuum DoS \eqref{eq:cont_Mathieu_DoS}. To make this precise, we consider the scaled discrete DoS
\begin{equation}
    \text{DoS}_{d}^{\epsilon}(\mu) := \frac{1}{\epsilon} \lim_{M \rightarrow \infty} \frac{1}{NM} \Tr \left( 1 + \left(\frac{H_d}{\epsilon^2} - \mu \right)^2 \right)^{-1} = \frac{1}{2\pi} \int_{-\pi}^{\pi}\,\Tr\left(1+\left(\frac{D_d(k)}{\epsilon^2} - \mu \right)^2\right)^{-1} \text{d}k.
\end{equation}

Our goal now is to prove the following theorem:
\begin{theorem} \label{thm:Mathieu_theorem}
For each fixed $\mu \in \mathbb{R}$, there exist positive constants $C > 0$ and $\epsilon_0 > 0$ such that, for all $0 < \epsilon < \epsilon_0$, we have
\begin{equation}
    \left| \text{\emph{DoS}}_d^\epsilon(\mu) - \text{\emph{DoS}}_c(\mu) \right| \leq C \epsilon^{\frac{1}{11}}.
\end{equation}
In particular,
\begin{equation}
    \lim_{\epsilon \rightarrow 0} \text{\emph{DoS}}_d^\epsilon(\mu) = \text{\emph{DoS}}_c(\mu).
\end{equation}
\end{theorem}

\begin{figure}[!htb]
\centering
\includegraphics[width=\columnwidth]{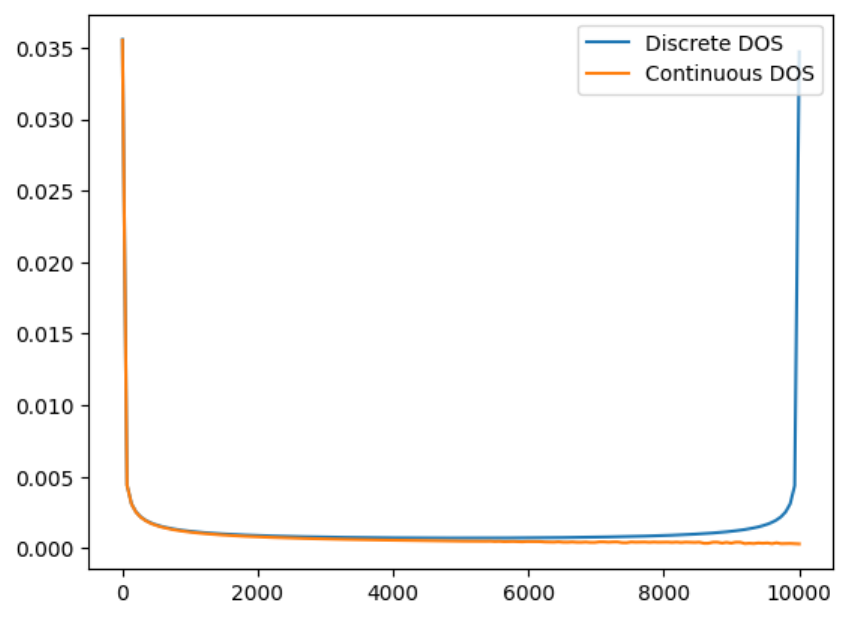}
    \caption{\label{fig:Mathieu Dos} Numerically computed re-scaled density of states for the discrete Mathieu equation \eqref{eq:disc_Mathieu_DoS} with $\lambda = 8$, compared with the numerically computed density of states for the continuum approximation \eqref{eq:cont_Mathieu_DoS}, showing the expected agreement near the bottom of the spectrum.}
\end{figure}

To see that this is plausible, note that the off diagonal terms of $D_c(k)$ and $\frac{D_d(k)}{\epsilon^2}$ are the same and the diagonal terms of $\frac{D_d(k)}{\epsilon^2}$ are $\frac{2 - 2 \cos(\epsilon (k+2\pi n))}{\epsilon^2}$ while the diagonal terms of $D_c(k)$ are $(k+2\pi n)^2$. Taylor expansion shows $\frac{2 - 2 \cos(\epsilon (k+2\pi n))}{\epsilon^2} \approx (k+2\pi n)^2 + O(\epsilon^2 k^4)$. Therefore in the limit as $\epsilon\downarrow0$ $D_c(k) \approx \frac{D_d(k)}{\epsilon^2}$.
However, this does not constitute a full proof. We provide a detailed proof in Appendix \ref{sec:mathieu_proofs}. Numerical results for the density of states calculated using both formulae are shown in Figure \ref{fig:Mathieu Dos}.

\newpage

\appendix

\section{Proof of Theorem \ref{thm:Laplace_theorem}} \label{sec:laplacian_proofs}

In this section we prove that 
\begin{equation}
\begin{split}
    \lim_{\epsilon \rightarrow 0} \text{DoS}_d^\epsilon(\mu) &= \lim_{\epsilon \rightarrow 0} \frac{1}{2\pi\epsilon}\int_{-\pi}^{\pi}\,{\frac{1}{1+\left(\frac{2 - 2\cos(k)}{\epsilon^2} - \mu\right)^2}\text{d}k} \\ &= \frac{1}{2\pi}\int_{-\infty}^{\infty}\,{\frac{1}{1+(k^2 - \mu)^2}\text{d}k} = \text{DoS}_c(\mu).
\end{split}
\end{equation}
for all $\mu \in \mathbb{R}$.

For any $q > 0$ the discrete integral can be split into three parts that can be evaluated separately
\begin{equation}
    \frac{1}{\epsilon}\int_{-\pi}^{\pi}\,{\frac{1}{1+\left(\frac{2 - 2\cos(k)}{\epsilon^2} - \mu\right)^2}\text{d}k} = I_1 + I_2 + I_3,
\end{equation}
where we split the integral into ``outer'' integrals
\begin{equation}
    I_1 := \frac{1}{\epsilon}\int_{-\pi}^{-\epsilon^{q/2}}\,{\frac{1}{1+\left(\frac{2 - 2\cos(k)}{\epsilon^2} - \mu\right)^2}\text{d}k}, \quad I_3 := \frac{1}{\epsilon}\int_{\epsilon^{q/2}}^{\pi}\,{\frac{1}{1+\left(\frac{2 - 2\cos(k)}{\epsilon^2} - \mu\right)^2}\text{d}k},
\end{equation}
and an ``inner'' integral
\begin{equation}
    I_2 := \frac{1}{\epsilon}\int_{-\epsilon^{q/2}}^{\epsilon^{q/2}}\,{\frac{1}{1+\left(\frac{2 - 2\cos(k)}{\epsilon^2} - \mu\right)^2}\text{d}k}.
\end{equation}

The outer integrals approach 0 as $\epsilon \downarrow 0$. Since the integrand is even, we show this WLOG for $I_3$ by showing that the largest value of the Lorentzian is at the bottom of the interval. From Taylor expansion we know 
\begin{equation}
    2 - 2\cos(k) = k^2 + O(k^4).
\end{equation}
For sufficiently small $k$ which includes $\epsilon^{\frac{q}{2}}$
\begin{equation}
    2 - 2\cos(k) \geq \frac{k^2}{2}.
\end{equation}
Since $\cos$ is monotone on the interval $[0,\pi]$, for all $k \in [\epsilon^{\frac{q}{2}},\pi]$
\begin{equation}
    2 - 2\cos(k) \geq 2 - \cos(\epsilon^{\frac{q}{2}}) \geq \frac{\epsilon^{q}}{2}.
\end{equation}
For $\mu<\frac{\epsilon^{q-2}}{4}$ we therefore have the bound
\begin{equation}
    \frac{2 - 2\cos(k)}{\epsilon^2} - \mu \geq \frac{\epsilon^{q-2}}{4}.
\end{equation}
This will hold for all $\mu$ in the limit if $q < 2$. This means the entire integrand is bounded

\begin{equation}
    \frac{1}{1+\left(\frac{2 - 2\cos(k)}{\epsilon^2} - \mu\right)^2} \leq 16\epsilon^{4-2q}.
\end{equation}
Therefore the entire integral is bounded
\begin{equation}
    \frac{1}{\epsilon}\int_{\epsilon^{q/2}}^{\pi}\,{\frac{1}{1+\left(\frac{2 - 2\cos(k)}{\epsilon^2} - \mu\right)^2}\text{d}k} \leq 16\pi \epsilon^{3-2q}.
\end{equation}
In the limit as $\epsilon \downarrow 0$ the integral is 0 assuming $q<\frac{3}{2}$. 

For the middle integral rescaling the bounds yields 
\begin{equation}
    \frac{1}{\epsilon}\int_{-\epsilon^{\frac{q}{2}}}^{\epsilon^{\frac{q}{2}}}\,{\frac{1}{1+\left(\frac{2 - 2\cos(k)}{\epsilon^2} - \mu\right)^2}\text{d}k} = \int_{-\epsilon^{\frac{q}{2}-1}}^{\epsilon^{\frac{q}{2}-1}}\,{\frac{1}{1+\left(\frac{2 - 2\cos(k)}{\epsilon^2} - \mu\right)^2}\text{d}k}.
\end{equation}
To compare to the continuum limit note that in the limit as $\epsilon \downarrow 0$
\begin{equation}
    \int_{\epsilon^{\frac{q}{2}-1}}^{\infty}\,{\frac{1}{1+(k^2 - \mu)^2}\text{d}k} = \int_{-\infty}^{-\epsilon^{\frac{q}{2}-1}}\,{\frac{1}{1+(k^2 - \mu)^2}\text{d}k} \rightarrow 0,
\end{equation}
since for $\mu \leq \frac{\epsilon^{\frac{q}{2}-1}}{2}$
\begin{equation}
    \int_{\epsilon^{\frac{q}{2}-1}}^{\infty}\,{\frac{1}{1+(k^2 - \mu)^2}\text{d}k} \leq \int_{\epsilon^{\frac{q}{2}-1}}^{\infty}\,{\frac{1}{(\frac{k^2}{2})^2}\text{d}k} = \frac{4\epsilon^{3-\frac{3q}{2}}}{3}\rightarrow 0.
\end{equation}
Therefore if 
\begin{equation}
    \int_{-\epsilon^{\frac{q}{2}-1}}^{\epsilon^{\frac{q}{2}-1}}\,{\frac{1}{1+\left(\frac{2 - 2\cos(k)}{\epsilon^2} - \mu\right)^2} - \frac{1}{1+(k^2 - \mu)^2}}\text{d}k
\end{equation}
approaches 0 the total difference between the discrete and continuum expression for DoS approaches 0. Rewriting the difference yields
\begin{equation}
    \int_{-\epsilon^{\frac{q}{2}-1}}^{\epsilon^{\frac{q}{2}-1}}\,{\frac{1}{1+\left(\frac{2 - 2\cos(k)}{\epsilon^2} - \mu\right)^2} - \frac{1}{1+(k^2 - \mu)^2}}\text{d}k = \int_{-\epsilon^{\frac{q}{2}-1}}^{\epsilon^{\frac{q}{2}-1}}\,{\frac{(k^2 - \mu)^2-\left(\frac{2 - 2\cos(k)}{\epsilon^2} - \mu\right)^2}{\left(1+\left(\frac{2 - 2\cos(k)}{\epsilon^2} - \mu\right)^2\right)(1+(k^2 - \mu)^2)}}\text{d}k.
\end{equation}
We note that the denominator is bounded below as
\begin{equation}
    \left(1+\left(\frac{2 - 2\cos(\epsilon k)}{\epsilon^2} - \mu\right)^2\right)(1+(k^2 - \mu)^2) \geq 1,
\end{equation}
and the numerator is bounded above as
\begin{equation}
\begin{split}
    (k^2 - \mu)^2-\left(\frac{2 - 2\cos(\epsilon k)}{\epsilon^2} - \mu\right)^2 &= k^4 - 4\epsilon^{-4}(1 - \cos(\epsilon k))^2 + 4\mu\epsilon^{-2}(1 - \cos(\epsilon k)) - 2\mu k^2 \\
    &\leq \epsilon^2k^6 + 2\epsilon^2k^4\mu \leq \epsilon^{3q-4} + 2\mu\epsilon^{2q-2}.
\end{split}
\end{equation}
Therefore
\begin{equation}
\begin{split}
    \int_{-\epsilon^{\frac{q}{2}-1}}^{\epsilon^{\frac{q}{2}-1}}\,{\frac{(k^2 - \mu)^2- \left( \frac{2 - 2\cos(k)}{\epsilon^2} - \mu \right)^2}{\left(1+\left(\frac{2 - 2\cos(k)}{\epsilon^2} - \mu \right)^2 \right)(1+(k^2 - \mu)^2)}}\text{d}k \leq \int_{-\epsilon^{\frac{q}{2}-1}}^{\epsilon^{\frac{q}{2}-1}}\,{\epsilon^{3q-4} + 2\mu \epsilon^{2q-2}}\text{d}k
    = 2\epsilon^{\frac{7q}{2}-5} + 4\mu\epsilon^{\frac{5q}{2}-3}.
\end{split}
\end{equation}
Incorporating the bounds on all portions of the integrals yields
\begin{equation} \label{eq:all_error_terms}
\begin{split}
    \frac{1}{\epsilon}\int_{-\pi}^{\pi}\,{\frac{1}{1+\left(\frac{2 - 2\cos(k)}{\epsilon^2} - \mu\right)^2}\text{d}k} - \int_{-\infty}^{\infty}\,{\frac{1}{1+(k^2 - \mu)^2}\text{d}k} \leq 2\epsilon^{\frac{7q}{2}-5} + 4\mu\epsilon^{\frac{5q}{2}-3} + 32\pi \epsilon^{3-2q} + \frac{8\epsilon^{3-\frac{3q}{2}}}{3}.
\end{split}
\end{equation}

For $\frac{10}{7}< q< \frac{3}{2}$ and $|\mu| \leq \epsilon^{q-2}$ this approaches 0 as $\epsilon \downarrow 0$. This again holds for all $\mu$ in the limit since $q < \frac{3}{2}$. The fastest decay occurs by balancing the first and third terms in \eqref{eq:all_error_terms} obtaining
\begin{equation}
    \frac{7}{2} q - 5 = 3 - 2 q \implies q = \frac{16}{11}. 
\end{equation}
These terms then contribute the dominant error term proportional to $\epsilon^{\frac{1}{11}}$. 

\section{Proof of Theorem \ref{thm:Mathieu_theorem}}
\label{sec:mathieu_proofs}
In this section we prove that
\begin{equation}
\begin{split}
    \lim_{\epsilon \rightarrow 0} \DoS_d^\epsilon(\mu) &= \lim_{\epsilon \rightarrow 0} \frac{1}{2 \pi} \int_{-\pi}^{\pi}\,\Tr\left(1+\left(\frac{D_d(k)}{\epsilon^2} - \mu \right)^2\right)^{-1} \text{d}k \\ &= \frac{1}{2\pi}\int_{-\pi}^{\pi}\,{\Tr((1+(D_c(k) - \mu )^2)^{-1} } )\text{d}k = \DoS_c(\mu),
\end{split}
\end{equation}
where $D_d(k)$ is defined by \eqref{eq:matheiu H_d(k)}, $D_c(k)$ by \eqref{eq:mathieu H_c(k)}, and $\mu \in \mathbb{R}$ is fixed but arbitrary.
To show that these integrals are equal we will prove that it is only necessary to compare the traces of smaller matrices than in the definitions. We will first prove that the trace in the continuous definition can be approximated by the trace of a smaller matrix.

\subsection{Continuous proof portion}

In this section we will prove the following lemma
\begin{lemma}\label{lemma: mathieu continuous}
Let $D_c^L(k)$ be defined as the following truncation of $D_c(k)$
\begin{equation} \label{eq:D_c' def}
    (D_c^L(k))_{n,m} := \begin{cases}
(D_c(k))_{n,m}, & \text{if } |n| < L \text{ and } |m| < L\\
0, & \text{otherwise}
\end{cases}
\end{equation}
with entries of the matrix measured from the middle using the following convention
\begin{equation}
    (D_c(k))_{n,n} := (k+2\pi n)^2.
\end{equation}
We can write the trace of an operator $A$ with respect to an orthonormal basis $\{ \ket{n} \}_{n \in \mathbb{Z}}$ as
\begin{equation}
   \Tr(A) := \sum_{n = -\infty}^{\infty}\,{\Braket{n|A|n}},
\end{equation}
whenever the sum converges. With $\{ \ket{n} \}_{n \in \mathbb{Z}}$ the Fourier basis as in \eqref{eq:matrix FH_dF_-1} we have
\begin{equation} \label{continuous lemma}
    \Tr((1+(D_c(k) - \mu )^2)^{-1}) = \lim_{L\rightarrow\infty} \sum_{n = -L}^{L}\,{\Braket{n|(1+(D_c^{2L}(k) - \mu)^2)^{-1}|n}}.
\end{equation}
In particular, the infinite sum defining the left-hand side exists, and can be computed by the limit on the right-hand side, which involves only the truncated operator \eqref{eq:D_c' def} with $L$ replaced by $2 L$.
\end{lemma}
To prove this it is necessary to break the Lorentzian of $D_c$ into pieces for which the traces can be calculated separately
\begin{equation}
   \Tr((1+(D_c(k) - \mu )^2)^{-1}) =\Tr((D_c(k)-\mu-i)^{-1}) -\Tr((D_c(k)-\mu+i)^{-1}),
\end{equation}
and then to break the calculation of the trace into a calculation of the traces of different portions of the matrix. 
We divide the trace into an upper, middle and lower portion 
\begin{equation}
\begin{split}
    \Tr((D_c(k)-\mu-i)^{-1}) = T_1 + T_2 + T_3
\end{split}
\end{equation}
where the upper portion of the trace is 
\begin{equation}
    T_1:= \sum_{n = -\infty}^{-L-1}\,{\Braket{n|(D_c(k)-\mu-i)^{-1}|n}}
\end{equation}
the middle portion of the trace is 
\begin{equation}
    T_2:= \sum_{n = -L}^{L}\,{\Braket{n|(D_c(k)-\mu-i)^{-1}|n}}
\end{equation}
and the lower portion of the trace is 
\begin{equation}
    T_3:= \sum_{n = L+1}^{\infty}\,{\Braket{n|(D_c(k)-\mu-i)^{-1}|n}}.
\end{equation}
We will show that the upper portion of the trace depends only on the upper portion of $(D_c(k)-\mu-i)^{-1}$ in the limit as $L\uparrow \infty$ and use that to show that the upper portion of the trace approaches 0 as $L\uparrow \infty$. We will do the same for the lower portion of the trace and we will similarly show that the middle portion of the trace only depends on the middle of $(D_c(k)-\mu-i)^{-1}$ in the same limit. 
We define the lower portion of the matrix as 
\begin{equation}
    (D_c'^L(k))_{n,m} := \begin{cases}
0, & \text{if } n < L \text{ or } m < L\\
(D_c(k))_{n,m}, & \text{otherwise}
\end{cases}.
\end{equation}
We can show that only this lower portion of the matrix contributes to the lower portion of the trace by showing the contribution of the upper portion of the matrix is 0. We rewrite the lower trace
\begin{equation}
    \begin{split}
        &\sum_{n = L+1}^{\infty}\,{\Braket{n|(D_c(k)-\mu-i)^{-1}|n}} = \\ 
        &\sum_{n = L+1}^{\infty}\,{\Braket{n|\left(D_c'^{\frac{L}{2}}(k)-\mu-i\right)^{-1}|n}} + \sum_{n = L + 1}^{\infty}\,{\Braket{n|(D_c(k)-\mu-i)^{-1} - \left(D_c'^{\frac{L}{2}}(k)-\mu-i\right)^{-1}|n}}.
    \end{split}
\end{equation}
We then rewrite the upper portion of the matrix
\begin{equation}
    (D_c(k)-\mu-i)^{-1} - \left(D_c'^{\frac{L}{2}}(k)-\mu-i\right)^{-1} = (D_c(k)-\mu-i)^{-1} \left(D_c(k)-D_c'^{\frac{L}{2}}(k)\right) \left(D_c'^{\frac{L}{2}}(k)-\mu-i\right)^{-1}.
\end{equation}
Writing out the matrix multiplication explicitly
\begin{equation} \label{continuous lower matrix multiplication}
    \begin{split}
        &\sum_{n = L + 1}^{\infty}\,{\Braket{n|(D_c(k)-\mu-i)^{-1} - \left(D_c'^{\frac{L}{2}}(k)-\mu-i\right)^{-1}|n}} = \\ &\sum_{n = L + 1}^{\infty}\, \sum_{m}(D_c(k)-\mu-i)^{-1}_{n,m}\sum_{m'}\left(D_c(k)-D_c'^{\frac{L}{2}}(k)\right)_{m,m'}\left(D_c'^{\frac{L}{2}}(k)-\mu-i\right)^{-1}_{m',n}.
    \end{split}
\end{equation}
We note
\begin{equation}
    \left(D_c(k)-D_c'^{\frac{L}{2}}(k)\right)_{m,m'} = 0 \text{ if } m> \frac{L}{2} \text{ and } m'> \frac{L}{2}.
\end{equation}
The off-diagonal terms of $D_c$ and $D_c'^L$ satisfy 
\begin{equation}
    D_c(k)_{n,m} \leq e\lambda e^{-|n-m|}, \quad \left(D_c'^{\frac{L}{2}}(k)\right)_{n,m} \leq e\lambda e^{-|n-m|}.
\end{equation}
Therefore by a Combes-Thomas-type estimate \cite{1973CombesThomas,doi:10.1137/15M1022628,WeinanLu2011} there exist constants $c_1,c_2 > 0$ independent of $L$ such that 
\begin{equation}
    (D_c(k)-\mu-i)^{-1}_{n,m} \leq c_1e^{-c_2|n-m|}, \quad \left(D_c'^{\frac{L}{2}}(k)-\mu-i\right)^{-1}_{n,m} \leq c_1e^{-c_2|n-m|}.
\end{equation}
Therefore for $m>\frac{L}{2}$ and $m'>\frac{L}{2}$
\begin{equation}
    (D_c(k)-\mu-i)^{-1}_{n,m}\left(D_c(k)-D_c'^{\frac{L}{2}}(k)\right)_{m,m'}\left(D_c'^{\frac{L}{2}}(k)-\mu-i\right)^{-1}_{m',n} = 0
\end{equation}
and for $m \leq \frac{L}{2}$ or $m' \leq \frac{L}{2}$
\begin{equation}
\begin{split}
    & (D_c(k)-\mu-i)^{-1}_{n,m}\left(D_c(k)-D_c'^{\frac{L}{2}}(k)\right)_{m,m'}\left(D_c'^{\frac{L}{2}}(k)-\mu-i\right)^{-1}_{m',n} 
    \\&\qquad \qquad \qquad \leq c_1^2 e^{-c_2|n-m|} e^{-c_2|n-m'|} \left(D_c(k)-D_c'^{\frac{L}{2}}(k)\right)_{m,m'}.
\end{split}
\end{equation}
Applying this to equation \eqref{continuous lower matrix multiplication} yields
\begin{equation}\label{explicit matrix multiplication lower continuous}
\begin{split}
    &\sum_{n = L + 1}^{\infty}\,{\Braket{n|(D_c(k)-\mu-i)^{-1} - \left(D_c'^{\frac{L}{2}}(k)-\mu-i\right)^{-1}|n}} \leq \\& \quad c_1^2 \sum_{n = L + 1}^{\infty}\, e^{-c_2n} \left( \sum_{m=-\infty}^{\infty}\sum_{m'=-\infty}^{\frac{L}{2}}e^{c_2m'}e^{-c_2\left|n-m\right|}\left(D_c(k)-D_c'^{\frac{L}{2}}(k)\right)_{m,m'} + \right. \\ & \quad \left. \sum_{m=-\infty}^{\frac{L}{2}}\sum_{m'= \frac{L}{2} +1 }^{\infty}e^{c_2m}e^{-c_2\left|n-m'\right|}\left(D_c(k)-D_c'^{\frac{L}{2}}(k)\right)_{m,m'} \right).
\end{split}
\end{equation}
For $|n-m|\geq 2$ we have $D_c(k)_{n,m} = 0$ and $D_c'^{\frac{L}{2}}(k)_{n,m} = 0$. Using this we have
\begin{equation}
\begin{split}
    = \; &c_1^2 \sum_{n = L + 1}^{\infty}\, e^{-2c_2n} \left( \sum_{m'=-\infty}^{\frac{L}{2}}\sum_{m = m' - 1}^{m'+1}e^{c_2(m+m')}\left(D_c(k)-D_c'^{\frac{L}{2}}(k)\right)_{m,m'} + \right. \\& \quad \left. e^{c_2(L+1)}\left(D_c(k)-D_c'^{\frac{L}{2}}(k)\right)_{\frac{L}{2},\frac{L}{2}+1}\right).
\end{split}
\end{equation}
Writing out the terms of the matrix we have 
\begin{equation}
    \sum_{m'=-\infty}^{\frac{L}{2}}\sum_{m=m'-1}^{m'+1}e^{c_2(m+m')}\left(D_c(k)-D_c'^{\frac{L}{2}}(k)\right)_{m,m'} = \left(\sum_{m=-\infty}^{\frac{L}{2}}\, e^{2c_2m}(k+2\pi m)^2 + 2\lambda e^{c_2(2m-1)}\right) + \lambda e^{c_2(L+1)},
\end{equation}
\begin{equation}
    e^{c_2(L+1)}\left(D_c(k)-D_c'^{\frac{L}{2}}(k)\right)_{\frac{L}{2},\frac{L}{2}+1} = \lambda e^{c_2(L+1)}.
\end{equation}
Therefore
\begin{equation}
\begin{split}
    \sum_{n = L + 1}^{\infty}\,{\Braket{n|(D_c(k)-\mu-i)^{-1} - \left(D_c'^{\frac{L}{2}}(k)-\mu-i\right)^{-1}|n}} \leq \\c_1^2 \sum_{n = L + 1}^{\infty}\,e^{-2c_2n}\sum_{m=-\infty}^{\frac{L}{2}}\, e^{2c_2m}(k+2\pi m)^2 + 2\lambda e^{c_2(2m+1)}.
\end{split}
\end{equation}
We can bound the sums over $n$ and $m$ separately. We first bound the sum over $m$. $\exists c_3>0$ independent of $L$ such that
\begin{equation}
    e^{-c_2m}(m+1)^2 \leq c_3 \text{ for } m>0.
\end{equation}
Therefore
\begin{equation}
    \sum_{m=-\infty}^{0}\, e^{2c_2m}(k+2\pi m)^2 \leq 4\pi^2 \sum_{m=-\infty}^{0}\, e^{c_2m} (e^{c_2m}(m-1)^2) \leq 4\pi^2 c_3\sum_{m= -\infty}^{0} e^{c_2m} = c_4 
\end{equation}
with $c_4$ also independent of $L$. For the upper part of the sum
\begin{equation}
    \sum_{m=1}^{\frac{L}{2}}\, e^{2c_2m}(k+2\pi m)^2 \leq 4\pi^2 \sum_{m=1}^{\frac{L}{2}}\, e^{2c_2m} \left(\frac{L}{2}+1\right)^2 = 4\pi^2\left(\frac{L}{2}+1\right)^2 \frac{e^{c_2L}-1}{e^{2c_2}-1} \leq 4 \pi^2 (e^{2c_2}-1)^{-1} L^2 e^{c_2L}.
\end{equation}
For the second term in the sum
\begin{equation}
  \sum_{m=-\infty}^{\frac{L}{2}}\, 2\lambda e^{c_2(2m+1)} = 2 \lambda \frac{e^{c_2(L+1)}}{e^{2c_2}-1} \leq 2 \lambda (e^{2c_2}-1)^{-1} e^{c_2L}.
\end{equation}
And for the sum over $n$
\begin{equation}
    \sum_{n = L + 1}^{\infty}\,e^{-2c_2n} = \frac{e^{-2c_2(L+1)}}{e^{2c_2} -1} 
\end{equation}
Substituting the bounds into equation \eqref{explicit matrix multiplication lower continuous} yields
\begin{multline}
    \sum_{n = L + 1}^{\infty}\,{\Braket{n|(D_c(k)-\mu-i)^{-1} - \left(D_c'^{\frac{L}{2}}(k)-\mu-i\right)^{-1}|n}} \leq \\ 2(2\pi^2L^2+\lambda) c_1^2 (e^{2c_2}-1)^{-2} e^{-c_2(L+2)} + c_4c_1^2(e^{2c_2}-1)^{-1} e^{-2c_2(L+1)}.
\end{multline}
Since $c_1$, $c_2$ and $c_4$ are independent of $L$ in the limit as $L$ approaches infinity
\begin{equation}
    \sum_{n = L + 1}^{\infty}\,{\Braket{n|(D_c(k)-\mu-i)^{-1} - \left(D_c'^{\frac{L}{2}}(k)-\mu-i\right)^{-1}|n}} = 0.
\end{equation}
To show that the lower portion of the trace is 0 we also need to show the lower portion of the matrix does not contribute to it. To do so we break that portion of the matrix into a diagonal and non-diagonal part. We define
\begin{equation}
    (B_c^L(k))_{n,m} := \begin{cases}
0, & \text{if } n \neq m\\
(D_c'^L(k))_{n,n} -\mu-i, & \text{otherwise}
\end{cases}
\end{equation}
\begin{equation}
    (A_c^L(k))_{n,m} := \begin{cases}
0, & \text{if } n = m\\
(D_c'^L(k))_{n,m}, & \text{otherwise}
\end{cases}.
\end{equation}
Writing the matrix in terms of these yields
\begin{equation*}
    D_c'^L(k) -\mu-i  = B_c^L(k) (I + (B_c^L(k))^{-1}A_c^L(k)).
\end{equation*}
We can express the inverse in terms of the diagonal and non-diagonal parts
\begin{equation*}
    (D_c'^L(k)-\mu-i)^{-1} = (B_c^L(k))^{-1} (I + (B_c^L(k))^{-1} A_c^L(k))^{-1}
\end{equation*}
\begin{equation*}
    ((B_c^L(k))^{-1})_{n,m} = \begin{cases}
0, & \text{if } n \neq m \text{ or } n<L\\
((k+2\pi n)^2 -\mu-i)^{-1}, & \text{otherwise}
\end{cases}.
\end{equation*}
Since we know the smallest eigenvalue of $B_c^L(k)$ and the largest eigenvalue of $A_c^L(k)$ we can show that the absolute value of all eigenvalues of $I + (B_c^L(k))^{-1} A_c^L(k)$ are greater than or equal to $1-\left|\frac{\lambda}{(k+2\pi L)^2 -\mu-i}\right|$. Therefore
\begin{equation*}
    \Braket{n|(I + (B_c^L(k))^{-1} A_c^L(k))^{-1}|n} \leq \left(1-\left|\frac{\lambda}{(k+2\pi L)^2 -\mu-i}\right|\right)^{-1} \leq 2.
\end{equation*}
Using this and knowing the eigenvalues of both $B_c^{\frac{L}{2}}(k)$ and $\left(B_c^{\frac{L}{2}}(k)\right)^{-1}$
\begin{equation}
\begin{split}
   \sum_{n = L + 1}^{\infty}\,{\Braket{n|\left(D_c'^{\frac{L}{2}}(k)-\mu-i\right)^{-1}|n}} = \sum_{n = L + 1}^{\infty}\,{\Braket{n|\left(B_c^{\frac{L}{2}}(k)\right)^{-1} \left(I + \left(B_c^{\frac{L}{2}}(k)\right)^{-1} A_c^{\frac{L}{2}}(k)\right)^{-1}|n}}
   \\ \leq \sum_{n = L + 1}^{\infty}\,{2|((k+2\pi n)^2 -\mu-i)^{-1}|} \leq 2 \sum_{n = L + 1}^{\infty}\,\frac{1}{n^2} \leq \frac{2}{L}.
\end{split}
\end{equation}
Therefore in the limit as $L$ approaches infinity
\begin{equation}
   \sum_{n = L + 1}^{\infty}\,{\Braket{n|\left(D_c'^{\frac{L}{2}}(k)-\mu-i\right)^{-1}|n}} \rightarrow 0.
\end{equation}
Since we have shown that in the limit both the upper and lower portions of the matrix do not contribute to the lower portion of the trace
\begin{equation}
   \sum_{n = L + 1}^{\infty}\,{\Braket{n|(D_c(k)-\mu-i)^{-1}|n}} \rightarrow 0
\end{equation}
and by the same logic 
\begin{equation}
    \sum_{n = -\infty}^{-L-1}\,{\Braket{n|(D_c(k)-\mu-i)^{-1}|n}} \rightarrow 0
\end{equation}
\begin{equation}
    \sum_{n = L + 1}^{\infty}\,{\Braket{n|(D_c(k)-\mu+i)^{-1}}} \rightarrow 0
\end{equation}
\begin{equation}
    \sum_{n = -\infty}^{-L-1}\,{\Braket{n|(D_c(k)-\mu+i)^{-1}}} \rightarrow 0.
\end{equation}
Therefore the infinite edges of the matrix have no contribution to the trace and as $L\uparrow\infty$
\begin{equation}
    \sum_{n = -L}^{L}\,{\Braket{n|(1+(D_c(k) - \mu )^2)^{-1}|n}} \rightarrow \Tr((1+(D_c(k) - \mu )^2)^{-1})
\end{equation}
exists.

Next we must show that only the middle of the matrix contributes to the middle portion of the trace. Recall the definition of the middle of the matrix from equation \ref{eq:D_c' def}
\begin{equation}
    (D_c^L(k))_{n,m} := \begin{cases}
(D_c(k))_{n,m}, & \text{if } |n| < L \text{ and } |m| < L\\
0, & \text{otherwise}
\end{cases}.
\end{equation}

We express the middle portion of the trace as the sum of the contribution of the middle of $D_c$ and the edges of $D_c$
\begin{multline}
    \sum_{n = -L}^{L}\,{\Braket{n|(D_c(k)-\mu-i)^{-1}|n}} = \sum_{n = -L}^{L}\,{\Braket{n|(D_c^{2L}(k) - \mu - i)^{-1}|n}} + \\ \sum_{n = -L}^{L}\,{\Braket{n|(D_c(k) - \mu -i)^{-1}-(D_c^{2L}(k) - \mu -i)^{-1}|n}}.
\end{multline}
To bound the contribution of the edges we can use the Combes Thomas estimate and that $D_c(k)=D_c^{2L}(k)$ for small $m$, $m'$ in the same way as before to rewrite the contribution of the edges as an explicit matrix multiplication. The result is
\begin{multline}\label{matrix multiplication middle continous}
    \sum_{n = -L}^{L}\,{\Braket{n|(D_c(k) - \mu -i)^{-1}-(D_c^{2L}(k) - \mu -i)^{-1}|n}} = \\ \sum_{n = -L}^{L}\, \sum_{m}(D_c(k)-\mu-i)^{-1}_{n,m}\sum_{m'}(D_c(k)-D_c^{2L}(k))_{m,m'}(D_c^{2L}(k)-\mu-i)^{-1}_{m',n} \\ \leq c_1^2 \sum_{n = -L}^{L}\,e^{2c_2n} 2 \left(\sum_{m=2L}^{\infty}\sum_{m'=m-1}^{m+1}e^{-c_2(m+m')}(D_c(k)-D_c^{2L}(k))_{m,m'} + e^{-c_2(4L-1)}\left(D_c(k)-D_c^{2L}(k)\right)_{2L-1,2L} \right).
\end{multline}
where $c_1$ and $c_2$ satisfy
\begin{equation}
    (D_c(k)-\mu-i)^{-1}_{n,m} \leq c_1e^{-c_2|n-m|}
\end{equation}
\begin{equation}
    \left(D_c^{2L}(k)-\mu-i\right)^{-1}_{n,m} \leq c_1e^{-c_2|n-m|}.
\end{equation}
Plugging in the entries of $D_c$ we can rewrite the sum over $m$
\begin{multline}
    \left(\sum_{m=2L}^{\infty}\sum_{m'=m-1}^{m+1}e^{-c_2(m+m')}(D_c(k)-D_c^{2L}(k))_{m,m'} + e^{-c_2(4L-1)}\left(D_c(k)-D_c^{2L}(k)\right)_{2L-1,2L} \right) = \\ \sum_{m=2L}^{\infty}e^{-2c_2m}(k+2\pi m)^2 + 2\lambda e^{-c_2(2m-1)}.
\end{multline}
We can bound the first term with
\begin{equation}
    \sum_{m=2L}^{\infty}e^{-2c_2m}(k+2\pi m)^2\leq c_3 \sum_{m=2L}^{\infty}e^{\frac{-3c_2m}{2}} = \frac{c_3 e^{-\frac{3c_2}{2}(2L-1)}}{e^{\frac{3c_2}{2}}-1}
\end{equation}
for some $c_3$ independent of $L$ satisfying
\begin{equation}
    e^{\frac{-c_2m}{2}}(k+2\pi m)^2 \leq c_3 \quad \forall m>0.
\end{equation}
We can explicitly evaluate the second term
\begin{equation}
  \sum_{m=2L}^{\infty}\, 2\lambda e^{-c_2(2m-1)} = 2 \lambda \frac{e^{-c_2(4L-1)}}{e^{2c_2}-1}. 
\end{equation}
We can do also bound the sum over $n$
\begin{equation}
    \sum_{n = -L}^{L}\,e^{2c_2n} = \frac{e^{2c_2L}-e^{-2c_2L}}{e^{2c_2}-1} \leq e^{2c_2L}(e^{2c_2}-1)^{-1}.
\end{equation}
Plugging both bounds back into \eqref{matrix multiplication middle continous} yields
\begin{multline}
    \sum_{n = -L}^{L}\,{\Braket{n|(D_c(k) - \mu -i)^{-1}-(D_c^{2L}(k) - \mu -i)^{-1}|n}} \leq \\ c_1^2(c_3(e^{\frac{3c_2}{2}}-1)^{-1}(e^{2c_2}-1)^{-1}e^{\frac{3c_2}{2}}e^{-c_2L} + 2\lambda(e^{2c_2}-1)^{-2}e^{c_2}e^{-2c_2L})
\end{multline}
which approaches $0$ as $L$ approaches $\infty$ since $c_1$, $c_2$ and $c_3$ are independent of $L$.
Similar logic shows 
\begin{equation}
    \sum_{n = -L}^{L}\,{\Braket{n|(D_c(k) - \mu +i)^{-1}-(D_c^{2L}(k) - \mu +i)^{-1}|n}} \rightarrow 0.
\end{equation}
From this and previous results showing the outer portion of the trace approaches 0 we have 
\begin{equation}
    \Tr((1+(D_c(k) - \mu )^2)^{-1}) = \lim_{L\rightarrow\infty} \sum_{n = -L}^{L}\,{\Braket{n|(1+(D_c^{2L}(k) - \mu)^2)^{-1}|n}}.
\end{equation}

Next we must prove that prove a similar lemma for $D_d(k)$.
\subsection{Discrete proof portion}

In this section we will prove the following lemma.
\begin{lemma}\label{lemma:mathieu discrete}
    Let $D_d^L$ be defined as the following truncation of $D_d$
\begin{equation} \label{eq:H_d' def}
    (D_d^L(k))_{n,m} := \begin{cases}
(D_d(k))_{n,m}, & \text{if } |n| < L \text{ and } |m| < L\\
0, & \text{otherwise}
\end{cases}.
\end{equation}
For $0<q<\frac{3}{2}$ we then have
\begin{equation} \label{discrete lemma}
    \lim_{\epsilon \rightarrow 0} \Tr\left(\left(1+\left(\frac{D_d(k)}{\epsilon^2} - \mu \right)^2\right)^{-1}\right) = \lim_{\epsilon \rightarrow 0+} \sum_{n = -\epsilon^{\frac{q}{2}-1}}^{\epsilon^{\frac{q}{2}-1}}\,{\Braket{n|\left( 1 + \left( \frac{D_d^{2\epsilon^{\frac{q}{2}-1}}(k)}{\epsilon^2} - \mu \right)^2 \right)^{-1}|n}}.
\end{equation}
\end{lemma}
Here we restrict $\epsilon$ and $q$ to values such that $\epsilon^{\frac{q}{2}-1}$ and $\frac{\epsilon^{-1}}{2}$ are integers so that we can use those values as matrix sizes and entries. 
Similarly to the continuous case we divide the matrix and the trace into segments. We define the lower portion of the matrix as 
\begin{equation}
    (D_d'^L(k))_{n,m} := \begin{cases}
0, & \text{if } n < L \text{ or } m < L\\
(D_d(k))_{n,m}, & \text{otherwise}
\end{cases}
\end{equation}
and divide the trace into 3 segments
\begin{equation}
    \begin{split}    \Tr\left(\left(1+\left(\frac{D_d(k)}{\epsilon^2} - \mu \right)^2\right)^{-1}\right) &= \sum_{n = \frac{-\epsilon^{-1}}{2}}^{-\epsilon^{\frac{q}{2}-1}-1}\,{\Braket{n|\left(1+\left(\frac{D_d(k)}{\epsilon^2} - \mu \right)^2\right)^{-1}|n}} \\&+ \sum_{n = -\epsilon^{\frac{q}{2}-1}}^{\epsilon^{\frac{q}{2}-1}}\,{\Braket{n|\left(1+\left(\frac{D_d(k)}{\epsilon^2} - \mu \right)^2\right)^{-1}|n}} \\ &+ \sum_{n = \epsilon^{\frac{q}{2}-1}+1}^{\frac{\epsilon^{-1}}{2}-1}\,{\Braket{n|\left(1+\left(\frac{D_d(k)}{\epsilon^2} - \mu \right)^2\right)^{-1}|n}}
\end{split}.
\end{equation}
We divide the lower segment into contributions from the upper and lower portions of $D_d$
\begin{multline}
   \sum_{n = \epsilon^{\frac{q}{2}-1}+1}^{\frac{\epsilon^{-1}}{2}-1}\,{\Braket{n|\left(1+\left(\frac{D_d(k)}{\epsilon^2} - \mu \right)^2\right)^{-1}|n}} = \sum_{n = \epsilon^{\frac{q}{2}-1}+1}^{\frac{\epsilon^{-1}}{2}-1}\,{\Braket{n|\left(1+\left(\frac{D_d'^{\frac{\epsilon^{\frac{q}{2}-1}}{2}}(k)}{\epsilon^2} - \mu \right)^2\right)^{-1}|n}} + \\ \sum_{n = \epsilon^{\frac{q}{2}-1}+1}^{\frac{\epsilon^{-1}}{2}-1}\,{\Braket{n|\left(\frac{D_d(k)}{\epsilon^2}-\mu-i\right)^{-1}-\left(\frac{D_d'^{\frac{\epsilon^{\frac{q}{2}-1}}{2}}(k)}{\epsilon^2} - \mu - i\right)^{-1}|n}} - \\ \sum_{n = \epsilon^{\frac{q}{2}-1}+1}^{\frac{\epsilon^{-1}}{2}-1}\,{\Braket{n|\left(\frac{D_d(k)}{\epsilon^2}-\mu+i\right)^{-1}-\left(\frac{D_d'^{\frac{\epsilon^{\frac{q}{2}-1}}{2}}(k)}{\epsilon^2} - \mu + i\right)^{-1}|n}}.
\end{multline}
We first show that the contribution from the lower portion of the matrix is 0. We do this by bounding the eigenvalues so we can bound the trace.
By the Gershgorin circle theorem \cite{varga2011gervsgorin}
\begin{equation*}
   1+\left(\frac{D_d'^{\frac{\epsilon^{\frac{q}{2}-1}}{2}}(k)}{\epsilon^2} - \mu \right)^2
\end{equation*}
has smallest eigenvalue at least
\begin{multline}
   1+ 2\lambda^2 + \left(2\epsilon^{-2}\left(1-\cos\left(\epsilon k + \pi \epsilon^{\frac{q}{2}}\right)\right)-\mu\right)^2 - 4\lambda\left(2\epsilon^{-2}\left(1-\cos\left(\epsilon k + \pi \epsilon^{\frac{q}{2}}\right)\right)-\mu\right) \\ \geq \frac{1}{2} \left(\left(k + 2\pi \epsilon^{\frac{q}{2}-1}\right)^2 -\mu\right)^2 \geq \frac{1}{8} \left(\pi \epsilon^{\frac{q}{2}-1}\right)^4.
\end{multline}
Therefore the largest eigenvalue of 
\begin{equation*}
   \left(1+\left(\frac{D_d'^{\frac{\epsilon^{\frac{q}{2}-1}}{2}}(k)}{\epsilon^2} - \mu \right)^2\right)^{-1}
\end{equation*}
is at most $\frac{8}{\pi^4} \epsilon^{4-2q}$.

This is also a bound for all the diagonal entries of the matrix. Therefore we can bound the lower trace by taking all the diagonal entries to be this value. Plugging this in as all of the diagonal entries yields
\begin{equation*}
   \sum_{n = \epsilon^{\frac{q}{2}-1}+1}^{\frac{\epsilon^{-1}}{2}-1}\,{\Braket{n|\left(1+\left(\frac{D_d'^{\frac{\epsilon^{\frac{q}{2}-1}}{2}}(k)}{\epsilon^2} - \mu \right)^2\right)^{-1}|n}} \leq \frac{\epsilon^{-1}}{2}\frac{8}{\pi^4} \epsilon^{4-2q}
\end{equation*}
which approaches 0 as $\epsilon$ approaches 0 for $q<\frac{3}{2}$.

For the contributions from the upper part of $D_d$ we write out the terms as an explicit matrix multiplication like in the continuous case
\begin{multline} \label{discrete lower matrix multiplication}
    \sum_{n = \epsilon^{\frac{q}{2}-1}+1}^{\frac{\epsilon^{-1}}{2}-1}\,{\Braket{n|\left(\frac{D_d(k)}{\epsilon^2}-\mu-i\right)^{-1}-\left(\frac{D_d'^{\frac{\epsilon^{\frac{q}{2}-1}}{2}}(k)}{\epsilon^2} - \mu - i\right)^{-1}|n}} = \\ \sum_{n = \epsilon^{\frac{q}{2}-1}+1}^{\frac{\epsilon^{-1}}{2}-1}\,{\Braket{n|\left(\frac{D_d(k)}{\epsilon^2}-\mu-i\right)^{-1}\frac{D_d(k) - D_d'^{\frac{\epsilon^{\frac{q}{2}-1}}{2}}(k)}{\epsilon^2}\left(\frac{D_d'^{\frac{\epsilon^{\frac{q}{2}-1}}{2}}(k)}{\epsilon^2} - \mu - i\right)^{-1}|n}} = \\ \sum_{n = \epsilon^{\frac{q}{2}-1}+1}^{\frac{\epsilon^{-1}}{2}-1}\,\sum_{m}\left(\frac{D_d(k)}{\epsilon^2}-\mu-i\right)^{-1}_{n,m}\sum_{m'}\left(\frac{D_d(k) - D_d'^{\frac{\epsilon^{\frac{q}{2}-1}}{2}}(k)}{\epsilon^2}\right)_{m,m'} \left( \frac{D_d'^{\frac{\epsilon^{\frac{q}{2}-1}}{2}}(k)}{\epsilon^2}-\mu-i\right)^{-1}_{m',n}.
\end{multline}
The same properties that allowed us to bound this contribution in the continuous case hold here
\begin{equation}
    \left(D_d(k) - D_d'^{\frac{\epsilon^{\frac{q}{2}-1}}{2}}(k)\right)_{m,m'} = 0 \text{ if } m> \frac{\epsilon^{\frac{q}{2}-1}}{2} \text{ or } m'> \frac{\epsilon^{\frac{q}{2}-1}}{2}.
\end{equation}
The off-diagonal terms of $D_d$ and $D_d'^L$ satisfy 
\begin{equation}
    D_d(k)_{n,m} \leq e\lambda \epsilon^2 e^{-|n-m|}
\end{equation}
\begin{equation}
    (D_d'^{L}(k))_{n,m} \leq e\lambda \epsilon^2 e^{-|n-m|}.
\end{equation}
Therefore by Combes-Thomas $\exists c_1,c_2>0$ independent of $\epsilon$ such that 
\begin{equation}
    \left(\frac{D_d(k)}{\epsilon^2}-\mu-i\right)^{-1}_{n,m} \leq c_1 e^{-c_2|n-m|}
\end{equation}
\begin{equation}
    \left(\frac{D_d'^{L}(k)}{\epsilon^2}-\mu-i\right)^{-1}_{n,m} \leq c_1 e^{-c_2|n-m|}.
\end{equation}
Therefore for $m>\frac{\epsilon^{\frac{q}{2}-1}}{2}$ and $m'>\frac{\epsilon^{\frac{q}{2}-1}}{2}$ or if $|m-m'|\geq 2$
\begin{equation}
    \Braket{n|\left(\frac{D_d(k)}{\epsilon^2}-\mu-i\right)^{-1}_{n,m}\left(\frac{D_d(k) - D_d'^{\frac{\epsilon^{\frac{q}{2}-1}}{2}}(k)}{\epsilon^2}\right)_{m,m'}\left(\frac{D_d'^{\frac{\epsilon^{\frac{q}{2}-1}}{2}}(k)}{\epsilon^2} - \mu - i\right)^{-1}_{m',n}|n} = 0,
\end{equation}
and for $m \leq \frac{\epsilon^{\frac{q}{2}-1}}{2}$ or $m' \leq \frac{\epsilon^{\frac{q}{2}-1}}{2}$
\begin{multline}
    \Braket{n|\left(\frac{D_d(k)}{\epsilon^2}-\mu-i\right)^{-1}_{n,m}\left(\frac{D_d(k) - D_d'^{\frac{\epsilon^{\frac{q}{2}-1}}{2}}(k)}{\epsilon^2}\right)_{m,m'}\left(\frac{D_d'^{\frac{\epsilon^{\frac{q}{2}-1}}{2}}(k)}{\epsilon^2} - \mu - i\right)^{-1}_{m',n}|n} \leq \\ c_1^2 e^{-c_2|n-m|} e^{-c_2|n-m'|} \left(\frac{D_d(k) - D_d'^{\frac{\epsilon^{\frac{q}{2}-1}}{2}}(k)}{\epsilon^2}\right)_{m,m'}.
\end{multline}
Using these bounds we can rewrite equation \eqref{discrete lower matrix multiplication}
\begin{multline}
    \sum_{n = \epsilon^{\frac{q}{2}-1}+1}^{\frac{\epsilon^{-1}}{2}-1}\,{\Braket{n|\left(\frac{D_d(k)}{\epsilon^2}-\mu-i\right)^{-1}-\left(\frac{D_d'^{\frac{\epsilon^{\frac{q}{2}-1}}{2}}(k)}{\epsilon^2} - \mu - i\right)^{-1}|n}} \leq \\ c_1^2 \sum_{n = \epsilon^{\frac{q}{2}-1}+1}^{\frac{\epsilon^{-1}}{2}-1}\, e^{-2c_2n} \left( \sum_{m=-\frac{\epsilon^{-1}}{2}+1}^{\frac{\epsilon^{\frac{q}{2}-1}}{2}}\sum_{m'=m-1}^{m+1}e^{c_2(m+m')}\left(\frac{D_d(k) - D_d'^{\frac{\epsilon^{\frac{q}{2}-1}}{2}}(k)}{\epsilon^2}\right)_{m,m'} + \right.\\ \left. e^{c_2\left(\epsilon^{\frac{q}{2}-1}+1\right)} \left(\frac{D_d(k) - D_d'^{\frac{\epsilon^{\frac{q}{2}-1}}{2}}(k)}{\epsilon^2}\right)_{\frac{\epsilon^{\frac{q}{2}-1}}{2}+1,\frac{\epsilon^{\frac{q}{2}-1}}{2}} + e^{c_2\epsilon^{-1}} \left(\frac{D_d(k) - D_d'^{\frac{\epsilon^{\frac{q}{2}-1}}{2}}(k)}{\epsilon^2}\right)_{-\frac{\epsilon^{-1}}{2},-\frac{\epsilon^{-1}}{2}} \right).
\end{multline}
Plugging in the entries of $D_d$ yields
\begin{multline}
    \sum_{m=-\frac{\epsilon^{-1}}{2}+1}^{\frac{\epsilon^{\frac{q}{2}-1}}{2}}\sum_{m'=m-1}^{m+1}e^{c_2(m+m')}\left(\frac{D_d(k) - D_d'^{\frac{\epsilon^{\frac{q}{2}-1}}{2}}(k)}{\epsilon^2}\right)_{m,m'} + \\ e^{c_2\left(\epsilon^{\frac{q}{2}-1}+1\right)} \left(\frac{D_d(k) - D_d'^{\frac{\epsilon^{\frac{q}{2}-1}}{2}}(k)}{\epsilon^2}\right)_{\frac{\epsilon^{\frac{q}{2}-1}}{2}+1,\frac{\epsilon^{\frac{q}{2}-1}}{2}} + e^{c_2\epsilon^{-1}} \left(\frac{D_d(k) - D_d'^{\frac{\epsilon^{\frac{q}{2}-1}}{2}}(k)}{\epsilon^2}\right)_{-\frac{\epsilon^{-1}}{2},-\frac{\epsilon^{-1}}{2}} \\= \sum_{m=-\frac{\epsilon^{-1}}{2}}^{\frac{\epsilon^{\frac{q}{2}-1}}{2}}\, e^{2c_2m}2\epsilon^{-2}(1-\cos(\epsilon k +2\pi m \epsilon)) + 2\lambda e^{c_2(2m+1)} \\ \leq \sum_{m=-\frac{\epsilon^{-1}}{2}}^{\frac{\epsilon^{\frac{q}{2}-1}}{2}}\, e^{2c_2m}8\epsilon^{-2} = 8 \epsilon^{-2} \frac{e^{-c_2\epsilon^{-1}}-e^{c_2\epsilon^{\frac{q}{2}-1}}}{1-e^{2c_2}} \leq 8 (e^{2c_2}-1)^{-1} \epsilon^{-2} e^{c_2\epsilon^{\frac{q}{2}-1}}.
\end{multline}
We now bound the sum over $n$
\begin{equation}
    \sum_{n = \epsilon^{\frac{q}{2}-1}+1}^{\frac{\epsilon^{-1}}{2}-1}\,e^{-2c_2n} = \frac{e^{-c_2\epsilon^{-1}}-e^{-2c_2\left(\epsilon^{\frac{q}{2}-1}+1\right)}}{e^{-2c_2} -1} \leq e^{-2c_2\left(\epsilon^{\frac{q}{2}-1}+1\right)}(e^{-2c_2} -1)^{-1}.
\end{equation}
Plugging in both bounds yields
\begin{equation}
    \sum_{n = \epsilon^{\frac{q}{2}-1}+1}^{\frac{\epsilon^{-1}}{2}-1}\,{\Braket{n|\left(\frac{D_d(k)}{\epsilon^2}-\mu-i\right)^{-1}-\left(\frac{D_d'^{\frac{\epsilon^{\frac{q}{2}-1}}{2}}(k)}{\epsilon^2} - \mu - i\right)^{-1}|n}} \leq 8 c_1^2 \epsilon^{-2} (e^{-c_2} - e^{c_2})^{-2} e^{-c_2\epsilon^{\frac{q}{2}-1}}.
\end{equation}
Since $c_1$ and $c_2$ are independent of $\epsilon$ in the limit as $\epsilon$ approaches 0 for $q<2$ 
\begin{equation}
    \sum_{n = \epsilon^{\frac{q}{2}-1}+1}^{\frac{\epsilon^{-1}}{2}-1}\,{\Braket{n|\left(\frac{D_d(k)}{\epsilon^2}-\mu-i\right)^{-1}-\left(\frac{D_d'^{\frac{\epsilon^{\frac{q}{2}-1}}{2}}(k)}{\epsilon^2} - \mu - i\right)^{-1}|n}} = 0.
\end{equation}
By the same logic
\begin{equation}
    \sum_{n = \epsilon^{\frac{q}{2}-1}+1}^{\frac{\epsilon^{-1}}{2}-1}\,{\Braket{n|\left(\frac{D_d(k)}{\epsilon^2}-\mu + i\right)^{-1}-\left(\frac{D_d'^{\frac{\epsilon^{\frac{q}{2}-1}}{2}}(k)}{\epsilon^2} - \mu + i\right)^{-1}|n}} \rightarrow 0.
\end{equation}
Since we have shown that the contribution to the lower trace from both the upper portion of $D_d$ and the lower portion are 0 for $q<\frac{3}{2}$
\begin{equation}
    \sum_{n = \epsilon^{\frac{q}{2}-1}+1}^{\frac{\epsilon^{-1}}{2}-1}\,{\Braket{n|\left(1+\left(\frac{D_d(k)}{\epsilon^2} - \mu \right)^2\right)^{-1}|n}} \rightarrow 0.
\end{equation}
By the same logic
\begin{equation}
    \sum_{n = \frac{-\epsilon^{-1}}{2}}^{-\epsilon^{\frac{q}{2}-1}-1}\,{\Braket{n|\left(1+\left(\frac{D_d(k)}{\epsilon^2} - \mu \right)^2\right)^{-1}|n}} \rightarrow 0.
\end{equation}
So only the middle segment of the trace is nonzero.

Next we must show that only the middle of $D_d$ contributes to the middle segment of the trace. Recall the definition of the middle of $D_d(k)$ from \ref{eq:H_d' def}
\begin{equation}
    (D_d^L(k))_{n,m} := \begin{cases}
(D_d(k))_{n,m}, & \text{if } |n| < L \text{ and } |m| < L\\
0, & \text{otherwise}
\end{cases}.
\end{equation}
We divide the middle trace into contributions from the middle and outer portions of $D_d$
\begin{multline}
    \sum_{n = -\epsilon^{\frac{q}{2}-1}}^{\epsilon^{\frac{q}{2}-1}}\,{\Braket{n|\left(\frac{D_d(k)}{\epsilon^2} -\mu-i \right)^{-1}|n}} = \sum_{n = -\epsilon^{\frac{q}{2}-1}}^{\epsilon^{\frac{q}{2}-1}}\,{\Braket{n|\left( \frac{D_d^{2\epsilon^{\frac{q}{2}-1}}(k)}{\epsilon^2} - \mu - i \right)^{-1}|n}} + \\ \sum_{n = -\epsilon^{\frac{q}{2}-1}}^{\epsilon^{\frac{q}{2}-1}}\,{\Braket{n|\left(\frac{D_d(k)}{\epsilon^2} -\mu-i \right)^{-1}-\left( \frac{D_d^{2\epsilon^{\frac{q}{2}-1}}(k)}{\epsilon^2} - \mu - i \right)^{-1}|n}}
\end{multline}
Using the same logic as in the continuous case it can be shown that as $\epsilon$ approaches 0 for q<2
\begin{equation}
    \sum_{n = -\epsilon^{\frac{q}{2}-1}}^{\epsilon^{\frac{q}{2}-1}}\,{\Braket{n|\left(\frac{D_d(k)}{\epsilon^2} -\mu-i \right)^{-1}-\left( \frac{D_d^{2\epsilon^{\frac{q}{2}-1}}(k)}{\epsilon^2} - \mu - i \right)^{-1}|n}} \rightarrow 0
\end{equation}
and similarly 
\begin{equation}
    \sum_{n = -\epsilon^{\frac{q}{2}-1}}^{\epsilon^{\frac{q}{2}-1}}\,{\Braket{n|\left(\frac{D_d(k)}{\epsilon^2} -\mu + i \right)^{-1}-\left( \frac{D_d^{2\epsilon^{\frac{q}{2}-1}}(k)}{\epsilon^2} - \mu + i \right)^{-1}|n}} \rightarrow 0.
\end{equation}
From this and earlier results we have 
\begin{equation}
    \Tr\left(\left(1+\left(\frac{D_d(k)}{\epsilon^2} - \mu \right)^2\right)^{-1}\right) = \lim_{\epsilon \rightarrow 0+} \sum_{n = -\epsilon^{\frac{q}{2}-1}}^{\epsilon^{\frac{q}{2}-1}}\,{\Braket{n|\left( 1 + \left( \frac{D_d^{2\epsilon^{\frac{q}{2}-1}}(k)}{\epsilon^2} - \mu \right)^2 \right)^{-1}|n}}
\end{equation}
for $0<q<\frac{3}{2}$.
Now that we can approximate the traces of the Lorentzians of $D_d$ and $D_c$ with portions of the traces of smaller matrices we must compare those matrices.
\subsection{Proof of difference portion}
From our earlier results we know that the only relevant part of the trace in both cases is the middle contribution to the middle segment of the trace. 
Therefore we now attempt to prove that the difference between those portions of the traces are 0. Specifically we will prove the following lemma.
\begin{lemma}\label{lemma:mathieu difference}
For $\frac{10}{7} < q$
\begin{equation}\label{difference lemma}
    \lim_{\epsilon \rightarrow 0+} \sum_{n = -\epsilon^{\frac{q}{2}-1}}^{\epsilon^{\frac{q}{2}-1}}\, \left( {\Braket{n|\left(1+\left(\frac{D_d^{2\epsilon^{\frac{q}{2}-1}}(k)}{\epsilon^2} - \mu \right)^2\right)^{-1}|n}} - {\Braket{n|\left(1+\left(D_c'^{2\epsilon^{\frac{q}{2}-1}}(k)-\mu\right)^2\right)^{-1}|n}} \right) = 0.
\end{equation}
\end{lemma}

To do this first note that 
\begin{equation}
    \begin{split} &\left(1+\left(\frac{D_d^{2\epsilon^{\frac{q}{2}-1}}(k)}{\epsilon^2} - \mu \right)^2\right)^{-1} - \left(1+\left(D_c'^{2\epsilon^{\frac{q}{2}-1}}(k)-\mu\right)^2\right)^{-1} = \left(1+\left(\frac{D_d^{2\epsilon^{\frac{q}{2}-1}}(k)}{\epsilon^2} - \mu \right)^2\right)^{-1}
    \\ &\left(\left(D_c'^{2\epsilon^{\frac{q}{2}-1}}(k)\right)^2 - \left(\frac{D_d^{2\epsilon^{\frac{q}{2}-1}}(k)}{\epsilon^2}\right)^2 + 2 D_c'^{2\epsilon^{\frac{q}{2}-1}}(k) \mu - 2 \frac{D_d^{2\epsilon^{\frac{q}{2}-1}}(k)}{\epsilon^2} \mu\right) \left(1+\left(D_c'^{2\epsilon^{\frac{q}{2}-1}}(k)-\mu\right)^2\right)^{-1},
    \end{split}
\end{equation}
and note that since
\begin{equation}
    1+\left(\frac{D_d^{2\epsilon^{\frac{q}{2}-1}}(k)}{\epsilon^2} - \mu \right)^2 
\end{equation}
and 
\begin{equation}
    1+\left(D_c'^{2\epsilon^{\frac{q}{2}-1}}(k)-\mu\right)^2
\end{equation}
both have all eigenvalues $\geq 1$ their inverses have all eigenvalues $\leq 1$.
Therefore the product
\begin{equation}
    \begin{split}    &\left(1+\left(\frac{D_d^{2\epsilon^{\frac{q}{2}-1}}(k)}{\epsilon^2} - \mu \right)^2\right)^{-1} *
    \\ &\left(\left(D_c'^{2\epsilon^{\frac{q}{2}-1}}(k)\right)^2 - \left(\frac{D_d^{2\epsilon^{\frac{q}{2}-1}}(k)}{\epsilon^2}\right)^2 + 2 D_c'^{2\epsilon^{\frac{q}{2}-1}}(k) \mu - 2 \frac{D_d^{2\epsilon^{\frac{q}{2}-1}}(k)}{\epsilon^2} \mu\right) *\\ &\left(1+\left(D_c'^{2\epsilon^{\frac{q}{2}-1}}(k)-\mu\right)^2\right)^{-1}.
    \end{split}
\end{equation}
has all diagonal entries less than or equal to the maximum eigenvalue of 
\begin{equation}
    \left(D_c'^{2\epsilon^{\frac{q}{2}-1}}(k)\right)^2 - \left(\frac{D_d^{2\epsilon^{\frac{q}{2}-1}}(k)}{\epsilon^2}\right)^2 + 2 D_c'^{2\epsilon^{\frac{q}{2}-1}}(k) \mu - 2 \frac{D_d^{2\epsilon^{\frac{q}{2}-1}}(k)}{\epsilon^2} \mu
\end{equation}
a matrix we can directly compute. 
Using the Gershgorin circle theorem \cite{varga2011gervsgorin} after calculating the terms of this matrix the maximum eigenvalue is less than
\begin{multline}
    \left(k+2\pi 2\epsilon^{\frac{q}{2}-1}\right)^4 - 4 \epsilon^{-4} \left(1-\cos\left(\epsilon\left(k+2\pi 2 \epsilon^{\frac{q}{2}-1}\right)\right)\right)^2  \\+ \mu\left(2 \left(k+2\pi 2\epsilon^{\frac{q}{2}-1}\right)^2 - 4 \epsilon^{-2} \left(1-\cos\left(\epsilon\left(k+2\pi 2 \epsilon^{\frac{q}{2}-1}\right)\right)\right)\right) \\ + 4\lambda \left(2 \left(k+2\pi 2\epsilon^{\frac{q}{2}-1}\right)^2 - 4 \epsilon^{-2} \left(1-\cos\left(\epsilon\left(k+2\pi 2 \epsilon^{\frac{q}{2}-1}\right)\right)\right)\right) \leq 2^{15} \pi^{6} \epsilon^{3q-4}.
\end{multline}
Therefore we can bound the sum of the diagonal entries
\begin{multline}
    \sum_{n = -\epsilon^{\frac{q}{2}-1}}^{\epsilon^{\frac{q}{2}-1}}\,{\Braket{n|\left(1+\left(\frac{D_d^{2\epsilon^{\frac{q}{2}-1}}(k)}{\epsilon^2} - \mu \right)^2\right)^{-1}|n}} - \sum_{n = -\epsilon^{\frac{q}{2}-1}}^{\epsilon^{\frac{q}{2}-1}}\,{\Braket{n|\left(1+\left(D_c'^{2\epsilon^{\frac{q}{2}-1}}(k)-\mu\right)^2\right)^{-1}|n}} \\ \leq 2 \epsilon^{\frac{q}{2}-1} 2^{15}\pi^{6} \epsilon^{3q-4} = 2^{16} \pi^{6} \epsilon^{\frac{7q}{2}-5}.
\end{multline}
This approaches 0 as $\epsilon$ approaches 0 if $q>\frac{10}{7}$, which is the lemma.

Now all that remains is to put these lemmas together.

\subsection{Proof of Theorem \ref{thm:Mathieu_theorem}}
From lemma \ref{lemma: mathieu continuous} for $q<2$ we have
\begin{equation}
    \lim_{\epsilon \rightarrow 0+} \sum_{n = -\epsilon^{\frac{q}{2}-1}}^{\epsilon^{\frac{q}{2}-1}}\, {\Braket{n|\left(1+\left(D_c'^{2\epsilon^{\frac{q}{2}-1}}(k)-\mu\right)^2\right)^{-1}|n}} = \Tr((1+(D_c(k) - \mu )^2)^{-1})
\end{equation}
since 
\begin{equation}
    \lim_{\epsilon \rightarrow 0+}f(\epsilon^{\frac{q}{2}-1}) = \lim_{N \rightarrow \infty} f(N).
\end{equation}
From lemma \ref{lemma:mathieu discrete} for $0<q<\frac{3}{2}$ we have 
\begin{equation}
     \lim_{\epsilon \rightarrow 0+} \sum_{n = -\epsilon^{\frac{q}{2}-1}}^{\epsilon^{\frac{q}{2}-1}}\,{\Braket{n|\left( 1 + \left( \frac{D_d^{2\epsilon^{\frac{q}{2}-1}}(k)}{\epsilon^2} - \mu \right)^2 \right)^{-1}|n}} = \lim_{\epsilon \rightarrow 0} \Tr\left(\left(1+\left(\frac{D_d(k)}{\epsilon^2} - \mu \right)^2\right)^{-1}\right).
\end{equation}
Therefore for $0<q<\frac{3}{2}$ we have 
\begin{equation}
\begin{split}
    & \lim_{\epsilon \rightarrow 0+} \sum_{n = -\epsilon^{\frac{q}{2}-1}}^{\epsilon^{\frac{q}{2}-1}}\, \left( {\Braket{n|\left(1+\left(\frac{D_d^{2\epsilon^{\frac{q}{2}-1}}(k)}{\epsilon^2} - \mu \right)^2\right)^{-1}|n}} - {\Braket{n|\left(1+\left(D_c'^{2\epsilon^{\frac{q}{2}-1}}(k)-\mu\right)^2\right)^{-1}|n}} \right) = \\ & \lim_{\epsilon \rightarrow 0} \Tr\left(\left(1+\left(\frac{D_d(k)}{\epsilon^2} - \mu \right)^2\right)^{-1}\right) - \Tr((1+(D_c(k) - \mu )^2)^{-1}).
\end{split}
\end{equation}
From lemma \ref{lemma:mathieu difference} for $q>\frac{10}{7}$ we have 
\begin{equation}
    \lim_{\epsilon \rightarrow 0+} \sum_{n = -\epsilon^{\frac{q}{2}-1}}^{\epsilon^{\frac{q}{2}-1}}\, \left( {\Braket{n|\left(1+\left(\frac{D_d^{2\epsilon^{\frac{q}{2}-1}}(k)}{\epsilon^2} - \mu \right)^2\right)^{-1}|n}} - {\Braket{n|\left(1+\left(D_c'^{2\epsilon^{\frac{q}{2}-1}}(k)-\mu\right)^2\right)^{-1}|n}} \right) = 0.
\end{equation}
Therefore for $\frac{10}{7}<q<\frac{3}{2}$ we have 
\begin{equation}
    \lim_{\epsilon \rightarrow 0} \Tr\left(\left(1+\left(\frac{D_d(k)}{\epsilon^2} - \mu \right)^2\right)^{-1}\right) - \Tr((1+(D_c(k) - \mu )^2)^{-1}) = 0
\end{equation}
and 
\begin{equation}
    \lim_{\epsilon \rightarrow 0} \text{DoS}_d^\epsilon(\mu) = \text{DoS}_c(\mu).
\end{equation}
This proves Theorem \ref{thm:Mathieu_theorem}. To determine a total bound on the convergence of $\DoS_c$ and $\DoS_d^{\epsilon}$ we combine all previously calculated bounds that decay slower than exponentially. For sufficiently small $\epsilon$ there are positive constants $c_5$ and $c_6$ such that 
\begin{equation}
    |\DoS_c(\mu)-\DoS_d^{\epsilon}(\mu)|\leq 2^{16} \pi^{6} \epsilon^{\frac{7q}{2}-5} + 8 \epsilon^{1-\frac{q}{2}} + \frac{8}{\pi^4}\epsilon^{3-2q} + c_5 \epsilon^{-2}e^{-c_6\epsilon^{\frac{q}{2}-1}}. 
\end{equation}
The fastest decay occurs for 
$q=\frac{16}{11}$ where the bound decays as 
$\epsilon^{\frac{1}{11}}$.

\bibliographystyle{plain}
\bibliography{ABW_library.bib}

\end{document}